\documentclass[preprint]{elsarticle}

\usepackage{epsfig,psfrag}
\usepackage{amsmath}
\usepackage{amssymb}
\usepackage{latexsym,pifont,color,comment}

\def\Bm{{\mathcal B}}

\def\A{\mathbb{A}}

\def\reals{{{\rm l} \kern -.15em {\rm R} }}

\def\A{{A}}
\def\Av{{\bf A}}
\def\B{{\bf B}}

\def\u{{\bf u}}
\def\E{{\bf E}}
\def\x{{\bf x}}
\def\n{{\bf n}}
\def\t{{\bf t}}
\def\r{{\bf r}}
\def\q{{q}}

\def\r{{r}}
\def\l{{\ell}}
\def\En{{\mathcal E}}
\def\half{\frac{1}{2}}
\def\Am{{\mathcal A}}
\def\Bm{{\mathcal B}}
\def\Cm{{\mathcal C}}

\def\be{\begin{equation}}
\def\ee{\end{equation}}

\def\ba{\begin{align}}
\def\ea{\end{align}}

\journal{Journal of Computational Physics}

\begin{document}

\begin{frontmatter}



\title{An Unstaggered Constrained Transport Method for the 3D Ideal Magnetohydrodynamic Equations}

\author[author1]{Christiane Helzel}
\ead{Christiane.Helzel@ruhr-uni-bochum.de}

\author[author2]{James A. Rossmanith\fnref{labc}}
\ead{rossmani@math.wisc.edu}

\author[author1]{Bertram Taetz}
\ead{Bertram.Taetz@rub.de}

\address[author1]{Fakult\"at f\"ur Mathematik, 
Ruhr-Universit\"at Bochum, 44780 Bochum, Germany}

\address[author2]{Department of Mathematics, University of Wisconsin,
480 Lincoln Drive, Madison, WI  53706-1388, USA}

\fntext[labc]{Corresponding author}

\begin{abstract}
Numerical methods for solving the ideal magnetohydrodynamic (MHD) 
equations in more than one space dimension must either confront the challenge
of controlling errors in the discrete divergence of the magnetic field, or else
be faced with nonlinear numerical instabilities. One approach for controlling
the discrete divergence is through a so-called {\it constrained transport}
method, which is based on first predicting a magnetic field through
a standard  finite volume solver, and then correcting
this field through the appropriate use of a magnetic vector potential.
In this work we
develop a constrained transport method for the 3D ideal MHD equations
that is based on a high-resolution wave propagation scheme. 
Our proposed scheme is the
3D extension of the 2D scheme developed by Rossmanith
[{\it SIAM J. Sci. Comp.} {\bf 28}, 1766 (2006)], and
is based on the high-resolution wave propagation method of
Langseth and LeVeque [{\it J. Comp. Phys.} {\bf 165}, 126 (2000)].
 In particular,
in our extension we take great care to maintain the three most important
properties of the 2D scheme: (1) all quantities, including all components of the
magnetic field and magnetic potential, are treated as cell-centered; (2)
we develop a high-resolution wave propagation scheme for evolving
the magnetic potential; and (3) we develop a wave limiting
approach that is applied during the vector potential evolution, which 
controls unphysical oscillations in the magnetic field.
One of the key numerical difficulties that is novel to 3D
is that the transport equation that must be solved
for the magnetic vector potential is only weakly hyperbolic.
In presenting our numerical algorithm we describe how to
numerically handle this problem of weak hyperbolicity, as
well as how to choose an appropriate gauge condition.
 The resulting scheme is applied to several
numerical test cases.
\end{abstract}

\end{frontmatter}
 
\section{Introduction}
\label{sec:intro}
The ideal magnetohydrodynamic (MHD) equations are a
common model for the macroscopic behavior of 
collisionless plasma  \cite{book:Chen84,book:Go98,book:Pa91}.
 These equations model the
fluid dynamics of an interacting mixture of positively and negatively 
charged particles, where each species is assumed to
behave as a charged ideal gas. Under the MHD assumption, this
mixture is taken to be quasi-neutral, meaning that as one
species moves, the other reacts instantaneously. This
 assumption allows one to collapse what should be
two sets of evolution equation
into a single set of equations for the total mass, momentum,
and energy of the mixture. Furthermore,  the resulting
dynamics are assumed to happen on slow time scales 
compared to the propagation time of light waves, which
yields a simplified set of Maxwell equations.
Finally, the term {\it ideal} refers to the fact that  
we assume the ideal Ohm's law: $\E = \B \times \u$,
which removes any explicit resistivity and the Hall
term\footnote{$\E$ is the electric field, $\B$ is the magnetic field,
and $\u$ is the macroscopic velocity of the plasma fluid.}.

All of the above described simplifications conspire to turn
the original two-fluid plasma model into a system that
can be viewed as a modified version of the compressible
Euler equations from gas dynamics. In particular, the ideal
MHD system can be written as a system of hyperbolic
conservation laws, where the conserved quantities 
are mass, momentum, energy, and magnetic field.
Furthermore, this system is equipped, just as the compressible
Euler equations are, with an entropy inequality that
features a convex scalar entropy and a corresponding
entropy flux. Indeed, the scalar entropy, with  some
help from the fact that the magnetic field is divergence-free, 
can be used to define entropy variables in which the 
MHD system is in symmetric hyperbolic form \cite{article:Barth05,article:Go72}.

As has been noted many times in the literature (e.g., 
Brackbill and Barnes \cite{article:BrBa80}, Evans and Hawley \cite{article:EvHa88},
and T\'oth \cite{article:To00}),
numerical methods for ideal MHD must in general satisfy (or at least control)
some discrete version of the divergence-free condition on the magnetic field:
\begin{equation}
\label{eqn:divBeq0}
\nabla \cdot \B = 0.
\end{equation}
Failure to accomplish this generically leads to a nonlinear numerical instability,
which often leads to negative pressures and/or densities.
Starting with the paper of Brackbill and Barnes \cite{article:BrBa80} in 1980,
several approaches for controlling errors in $\nabla \cdot \B$ have
been proposed. An in-depth review of many of these methods can
be found in T\'oth \cite{article:To00}. We very briefly summarize
the main points below.
\begin{description}
\item[{\bf Projection methods.}] 
Projection methods for ideal MHD are based on a predictor-corrector approach
for the magnetic field. Some standard finite volume or finite difference method is used 
to solve the ideal MHD
equations from time $t^n$ to $t^{n+1}$.
The approximate magnetic field that is predicted at time $t^{n+1}$ is denoted
by $\B^{\star}$. In general, $\B^{\star}$ contains both a nontrivial divergence-free
subspace and a divergent subspace; one would like to extract only the divergence-free
part and discard the divergent part. This can be accomplished by setting
\be
\label{eqn:project_new}
	\B^{n+1} := \B^{\star} - \nabla \psi,
\ee
and forcing $\nabla \cdot \B^{n+1} = 0$, which results
in a Poisson equation for the scalar potential $\psi$:
\be
	\nabla^2 \psi = \nabla \cdot \B^{\star}.
\ee
Once $\psi$ is computed, the divergence-free magnetic field at time $t^{n+1}$ is 
taken to be \eqref{eqn:project_new} (see T\'oth \cite{article:To00} and Balsara and Kim \cite{article:BaKi04} for further discussion).

This method is attractive for it allows a variety of methods to be
used in the prediction step, and then only requires one Poisson
solve per time-step to correct it. The clear disadvantage of this
approach is that it requires a global elliptic solve on a problem,
ideal MHD, that is purely hyperbolic. This could  be especially
computationally inefficient in the case of adaptively refined grids.

\item[{\bf The 8-wave formulation.}]  The linearized ideal MHD equations
support seven propagating plane-wave solutions\footnote{The seven propagating
waves are made up of two fast magnetosonic, two Alfv$\acute{\text{e}}$n, two slow magnetosonic,
and an entropy wave.} and a stationary plane-wave solution\footnote{
This is the so-called {\it divergence wave}.}.
The stationary plane-wave solution comes directly from the fact that
the (nonlinear) MHD equations  preserve the divergence constraint:
\begin{equation}
	\left( \nabla \cdot \B \right)_{,t} = 0,
\end{equation}
where $,t$ denotes the partial derivative with respect to time.
A seemingly unrelated fact was proved by Godunov \cite{article:Go72}:
the ideal MHD can only be put in symmetric hyperbolic form
if one adds to the ideal MHD equations a term that is 
proportional to $\nabla \cdot \B$ (in effect adding a term that,
at least on the continuous level, is  zero).
As it turns out, this additional ``source term'' not only allows
the equations to be put in symmetric hyperbolic form, but
it also restores Galilean invariance. With the additional term 
the divergence wave is no longer stationary,
and instead, it propagates with the fluid velocity:
\begin{equation}
 \left( \nabla \cdot \B \right)_{,t} +   \nabla \cdot \left( \u \nabla \cdot \B \right) = 0.
\end{equation}
From the point of view of numerical methods, the difference between a stationary and a
propagating divergence wave turns out to be very significant. Powell
\cite{article:Po94} and Powell et al. \cite{article:Po99} showed that numerical methods applied to the MHD
equations with the symmetrizing ``source term'' were 
much more stable than the same methods applied
to the original MHD equations.
The modified form of the MHD equations has come be called
the {\it 8-wave formulation}, since this form of the equations 
supports eight propagating plane wave solutions. 
Although this approach has been used with
some success (see Powell et al. \cite{article:Po99}), it does have a significant 
drawback: the 8-wave formulation is non-conservative and
difficulties with obtaining the correct weak solution have
been documented in the literature (see for example T\'oth \cite{article:To00}).

\item[{\bf Hyperbolic divergence cleaning methods.}] This method
was introduced by Dedner et al. \cite{article:Ded01a}, and is a 
close cousin to the above described projection method. 
The basic idea is to again solve for the divergence error
in the magnetic field. Instead of solving an elliptic equation,
however, a damped hyperbolic equation is prescribed for
the divergence error. This does not produce an exact divergence-free
magnetic field; however, it allows for the divergence error to be propagated
and damped away from where it originated. The main advantage of
this method is that it is easy to implement and requires no
elliptic solve. The main disadvantage is that this approach has
two tunable parameters: the speed of propagation of the error and
the rate at which the divergence error is damped. 

\item[{\bf Constrained transport methods.}]
The constrained transport (CT) approach for ideal MHD was introduced by Evans and Hawley
\cite{article:EvHa88}. The method is a modification of Yee's method \cite{article:Ye66} for electromagnetic 
wave propagation, and, at least in its original formulation, introduced staggered magnetic and electric fields.
Since the introduction of the CT framework, several variants and modifications 
have been introduced, including the work of Balsara \cite{article:Ba04},  Balsara and Spicer \cite{article:BaSp99a}, Dai and Woodward \cite{article:DaWo98}, 
Fey and Torrilhon \cite{article:FeTo03}, Londrillo and Del Zanna \cite{article:LoZa04},
Rossmanith \cite{article:Ro04b}, Ryu et al. \cite{article:Ryu98}, and De Sterck \cite{article:De01b}. An overview of many of these approaches, as well as the
introduction of a few more variants, can be found in
T\'oth \cite{article:To00}. In particular, T\'oth \cite{article:To00} showed
that a staggered magnetic field is not necessary. 

The CT framework, at least several versions of it, can also be viewed as a kind of
predictor-corrector approach for the magnetic field. Roughly speaking, the idea is
to again compute all of the conserved quantities with a standard finite volume method. From these computed quantities one
then constructs an approximation to the electric field through the ideal Ohm's law.
This electric field can then be used to update the magnetic vector potential,
which in turn, can be used to compute a divergence-free magnetic field (see
\S\ref{sec:framework} for more details).

The main advantages of this approach are that (1) there is no elliptic solve
and (2) there are no free parameters to choose such as in the hyperbolic
divergence-cleaning technique.
\end{description}

We focus in this work on 
developing a constrained transport method for the 3D ideal MHD equations
based on the high-resolution {\it wave propagation scheme} of 
LeVeque \cite{article:Le97} and its 3D extension 
developed by Langseth and LeVeque \cite{article:LaLe00}.
Our proposed scheme is the
3D extension of the 2D constrained transport scheme developed by Rossmanith
\cite{article:Ro04b}. 

The wave propagation scheme is an  unsplit finite volume method
that achieves second-order accuracy and optimal stability\footnote{The
wave propagation scheme is optimally stable in the sense that
no other fully explicit method that uses at most a 3-point stencil in 1D,
a 9-point stencil in 2D, and a 27-point stencil in 3D can achieve
larger maximum Courant numbers.}
in higher-dimensions through
the use of so-called {\it transverse} Riemann solvers. 
These transverse Riemann solvers are based partly on the work of
Colella \cite{article:Col90} on {\it corner transport upwind} (CTU) methods.
It is worth noting that there has been recent work on CTU methods
in the context of MHD by Gardiner and Stone
 \cite{article:athena05, article:athena08}
and Mignone and Tzeferacos \cite{article:MiTz10}.
Gardiner and Stone
 \cite{article:athena05, article:athena08} developed a constrained transport
 approach, while Mignone and Tzeferacos \cite{article:MiTz10} considered
 the hyperbolic divergence cleaning method of \cite{article:Ded01a} in the 
 context of the CTU scheme.
The method we propose in this work is therefore in the same class of
methods as that of  Gardiner and Stone \cite{article:athena05, article:athena08}
(i.e., unsplit finite volume methods with constrained transport),
and to a lesser extent Mignone and Tzeferacos \cite{article:MiTz10}
(i.e., unsplit methods for MHD), although
the details of our base scheme and especially the details of how
the magnetic field is updated differ greatly from these approaches.

In particular, the approach for updating the magnetic field that we propose in this work is
based on a generalization of the 2D constrained transport scheme of \cite{article:Ro04b},
which is equipped with all of the following features:
\begin{enumerate}
\item All quantities, including all components of the
magnetic field and magnetic potential, are treated as cell-centered;
\item The magnetic potential is evolved alongside the eight conserved
MHD variables via a modified non-conservative high-resolution wave propagation scheme; and 
\item Special limiters are applied in the evolution of the magnetic potential, which
control unphysical oscillations in the magnetic field.
\end{enumerate}
The scheme developed in this work extends all three of the above 
features to the 3D case. The new challenge that arises in 3D is that
the magnetic potential is a vector potential (instead of a scalar as in the 2D case).
Furthermore, this vector potential obeys, at least under gauge condition we advocate
in this work, a non-conservative {\it weakly hyperbolic} evolution equation.
We are thus faced with two important challenges:
\begin{enumerate}
\item Developing a  wave propagation scheme for	
	this non-conservative weakly hyperbolic evolution equation; and
\item Constructing appropriate limiters that again have the effect of
	limiting the magnetic potential in such a way as to
	produce a non-oscillatory magnetic field.
\end{enumerate}

After briefly reviewing the MHD equations in \S \ref{sec:equations},
we describe an overview of our proposed method in \S \ref{sec:framework}.
One important issue that arises from this approach is the gauge choice;
we discuss several possible choices in \S \ref{sec:vecpotential}.
Under the gauge condition that we choose, the 3D transport equation
that must be solved for the magnetic potential is only weakly hyperbolic.
We describe in detail in \S \ref{sec:numerics} how to numerically handle this difficulty.
 The resulting scheme is applied to several numerical test cases in \S \ref{sec:numex}.

\section{Basic equations}
\label{sec:equations}
The ideal magnetohydrodynamic (MHD) equations are a classical model
from plasma physics that describe the macroscopic evolution of a 
quasi-neutral two-fluid plasma system. Under the quasi-neutral
assumption, the two-fluid equations can be collapsed
into a single set of fluid equations for the total mass, momentum, and energy of the
system. The resulting equations can be written in the following 
form\footnote{We use boldface letters to denote vectors in physical space (i.e.,
$\reals^3$),
and $\| \cdot \|$ to denote the Euclidean norm of vector in the physical space.
Vectors in solution space, such as $q\in \reals^8$, 
where $q$ is the vector of conserved variables for the ideal MHD equations:
$q=\left(\rho,\, \rho \u, \, {\mathcal E}, \B \right)^T$, are not denoted using
boldface letters.}:
\begin{gather}
\label{eqn:MHD}
 \frac{\partial}{\partial t}
  \begin{bmatrix}
    \rho \\ \rho \u \\ \En \\ \B
  \end{bmatrix} +  \nabla \cdot
  \begin{bmatrix} \rho \u \\ \rho \u  \u + \left( {p} + \frac{1}{2}
  \| \B \|^2  \right) {\mathbb I}
    - \B  \B \\
    \u \left(\En + {p} + \frac{1}{2} \| \B \|^2 \right) - 
    	\B \left(\u \cdot \B \right)
    \\ \u  \B - \B  \u
  \end{bmatrix} = 0, \\
  \label{eqn:divfree}
  \nabla \cdot \B = 0,
\end{gather}
where $\rho$, $\rho \u$, and $\En$ are the total mass,  momentum,
and energy densities of the plasma system, and $\B$ is the
magnetic field. The thermal pressure, $p$, is related to the
conserved quantities through the ideal gas law:
\begin{equation}
\label{eqn:eos}
      p = \left( \gamma-1 \right) \left( \En - \frac{1}{2} \| \B \|^2 
      	- \frac{1}{2} \rho \| \u \|^2 \right),
\end{equation}
where $\gamma = 5/3$ is the ideal gas constant.

The equation for the magnetic field comes from
Faraday's law:
\begin{equation}
   \B_{,t} + \nabla \times \E = 0,
\end{equation}
where the electric field, $\E$, is approximated
by Ohm's law for a perfect conductor:
\begin{equation}
   \E = \B \times \u.
\end{equation}
Note that we have used in the above expression a comma followed by a subscript as  short-hand notation
for partial differentiation.  This notation is standard in many areas of mathematics,
 most notably in relativity theory (e.g., \cite{book:Poisson04}). 
We will continue to use this notation throughout this paper.

Under the Ohm's law assumption, we can
rewrite Faraday's law in the following
divergence form:
\begin{equation}
   \B_{,t} + \nabla \times \left( \B \times \u \right) = \B_{,t} + \nabla \cdot \left( \u \B - \B \u \right) = 0.
\end{equation}
Since the electric field is determined entirely from Ohm's law, we 
do not require an evolution equation for it; and thus, the only
other piece that we need from Maxwell's equations is
the divergence-free condition on the magnetic field
\eqref{eqn:divfree}. A complete derivation and discussion
of MHD system \eqref{eqn:MHD}-\eqref{eqn:divfree} can
be found in several standard plasma physics textbooks
(e.g., \cite{book:Chen84,book:Go98,book:Pa91}).

\subsection{$\nabla \cdot \B = 0$ is an involution}
We first note that system \eqref{eqn:MHD}, along with the
equation of state \eqref{eqn:eos}, provides
a full set of equations for the time evolution 
of all eight state variables: $\left( \rho, \, \rho \u, \, \En, \, \B \right)$.
These evolution equations form a hyperbolic system. In particular, the
eigenvalues of the flux Jacobian in some arbitrary
direction ${\bf n}$ ($\| {\bf n} \| = 1$) can be written as follows:
  \begin{alignat}{2}
  \lambda^{1,8} &= {\bf u} \cdot {\bf n} \mp c_f & &\text{ :  fast magnetosonic waves,} \\
  \lambda^{2,7} &= {\bf u} \cdot {\bf n}  \mp c_a & &\text{ :  Alfv$\acute{\text{e}}$n waves,} \\
  \lambda^{3,6} &= {\bf u} \cdot {\bf n}  \mp c_s & &\text{ :  slow magnetosonic waves,} \\
  \lambda^{4} &= {\bf u} \cdot {\bf n} & &\text{ :  entropy wave,} \\
  \lambda^{5} &= {\bf u} \cdot {\bf n}  & &\text{ :  divergence wave,}
\end{alignat}
where 
\begin{align}
  a &\equiv \sqrt{\frac{\gamma p}{\rho}}, \\
  c_a &\equiv \sqrt{\frac{\left(\B \cdot \n \right)^2}{\rho}}, \\
  c_f &\equiv \left\{ \frac{1}{2} \left[ a^2 +
      \frac{\|\B\|^{2}}{\rho} + \sqrt{\left(a^2 +
      \frac{\|\B\|^{2}}{\rho} \right)^{2} - 4 a^2
      \frac{\left(\B \cdot \n \right)^2}{\rho}} \right] \right\}^{1/2}, \\
  c_s &\equiv \left\{ \frac{1}{2} \left[ a^2 +
      \frac{\|\B\|^{2}}{\rho} - \sqrt{\left(a^2 +
      \frac{\|\B\|^{2}}{\rho} \right)^{2} - 4 a^2
      \frac{\left(\B \cdot \n \right)^2}{\rho}} \right] \right\}^{1/2} \, .
\end{align}
The eigenvalues are well-ordered in the sense that
\begin{equation}
  \lambda^{1} \le \lambda^{2} \le \lambda^{3} \le \lambda^{4} \le \lambda^{5} \le \lambda^{6}
  \le \lambda^{7} \le \lambda^{8} \, .
\end{equation}
The fast and slow magnetosonic waves are genuinely nonlinear, while the remaining waves are linearly degenerate.
Note that the so-called {\it divergence-wave} has been
made to travel at the speed $\u \cdot {\bf n}$ via the 
8-wave formulation described in \S\ref{sec:intro},
thus restoring Galilean invariance \cite{article:Go72,article:Po94,article:Po99}.
Note that despite the fact that we use the eigenvalues and eigenvectors
of the 8-wave formulation of the MHD equations, we will still
solve the MHD equations in conservative form (i.e., without the 
Godunov-Powell ``source term'').

The additional equation \eqref{eqn:divfree} is not needed 
in the time evolution of the conserved variables in the following sense:
\begin{quote}
	{\it If \eqref{eqn:divfree} is true at some time $t=T$, then
	evolution equation  \eqref{eqn:MHD} guarantees that
	 \eqref{eqn:divfree} is true for all time.}
\end{quote}
This result follows from taking the divergence of Faraday's law, which yields:
\begin{equation}
 \left(  \nabla \cdot \B \right)_{,t}  = 0.
\end{equation}
For this reason,  \eqref{eqn:divfree} should not be regarded as {\it constraint}
(such as the $\nabla \cdot \u =0$ constraint for the incompressible
Navier-Stokes equations), but rather an {\it involution} \cite{book:Dafermos10}.

\subsection{The role of $\nabla \cdot \B=0$ in numerical discretizations}
Although $\nabla \cdot \B = 0$ is an involution, and therefore has no
dynamic impact on the evolution of the exact MHD system, the story
is more complicated for numerical discretizations of ideal MHD.
Brackbill and Barnes \cite{article:BrBa80} gave a physical explanation
as to why $\nabla \cdot \B = 0$ should be satisfied in some appropriate
discrete sense: 
\begin{quote}
{\it   If $\nabla \cdot \B \ne 0$, then the magnetic force,
\be
	{\bf F} = \nabla \cdot \left\{ \B \B - \frac{1}{2} \| \B \|^2 \, {\mathbb I}
		\right\},
\ee
 in the direction of the
   magnetic field, will not in general vanish:
   \be  {\bf F} \cdot \B 
   	= \| \B \|^2 \, \left( \nabla \cdot \B\right) \ne 0. \ee
  If this spurious forcing becomes too large, it can lead to numerical instabilities (see
  for example \cite{article:BrBa80,article:Ro04b,article:To00}).}
\end{quote}

T\'oth \cite{article:To02} analyzed the magnetic force for central
finite difference methods and develop an approach
with the property that ${\bf F} \cdot \B=0$ on the discrete level.
This analysis makes it clear that simply having $\nabla \cdot \B = 0$ on the
discrete level is not sufficient to guarantee that ${\bf F} \cdot \B=0$ is
satisfied. Londrillo and Del Zanna \cite{article:LoZa04} came to the same conclusion
and proposed a modified approach for computing numerical fluxes. 
However, in \cite{article:Ro04b} it is argued that by satisfying an appropriate
discrete form of $\nabla \cdot \B=0$ the error in ${\bf F} \cdot \B$
is sufficiently controlled; and therefore,
the additional modifications
that are proposed in \cite{article:LoZa04} and \cite{article:To02} are not found to be
necessary. We omit the details here; and instead, refer the reader to 
\cite{article:LoZa04,article:Ro04b,article:To02}.

Another explanation as to why $\nabla \cdot \B=0$ should not be
ignored in numerical discretizations of MHD from a slightly different point-of-view, was offered
by Barth \cite{article:Barth05}.
Barth's explanation is based on the well-known result of Godunov \cite{article:Go72}
that the MHD entropy density,
\be
U(q) = -\rho \log \left( p \rho^{-\gamma} \right),
\ee
produces a set of {\it entropy variables}, $U_{,q}$, that do not immediately
symmetrize the ideal MHD equations. Instead, a symmetric hyperbolic
form of ideal MHD can only be obtained if an additional term
that is proportional to the divergence of the magnetic field
is included in the MHD equations:
\be
\underset{\text{ideal MHD}}{\underbrace{q_{,t} + {\nabla} \cdot {\bf F}(q)}} +  
		\underset{\text{additional term}}{\underbrace{\chi_{,q} \,
			{\nabla} \cdot \B}} = 0, \qquad \text{where} \quad
	\chi(q) = \left( \gamma -1 \right) \frac{\rho \u \cdot \B}{p}.
\ee
By looking at how the entropy behaves on the discrete
level, Barth \cite{article:Barth05} was able to prove that certain
discontinuous Galerkin discretizations of the ideal MHD
equtions could be made to be {\it entropy stable} (see
Tadmor \cite{article:Tad03}) if the discrete magnetic field where made globally
divergence-free.
The implication of this result is that schemes that do not
control errors in the divergence of the magnetic field
run the risk of becoming entropy unstable.

\section{An unstaggered constrained transport framework}
\label{sec:framework}
The constrained transport framework advocated in this
work is based on using the {\it magnetic vector potential},
$\Av : \reals^+ \times \reals^{3} \rightarrow \reals^3$, whose
curl yields the magnetic field:
\be
	\B = \nabla \times \Av.
\ee
A key feature of the proposed method is that 
the magnetic vector potential is updated alongside
the conserved variables: $q=(\rho, \rho \u, \En, \B)^T$.

The use of the magnetic potential in computing solutions
to the MHD equations goes all the way back to
Wilson \cite{article:Wilson75} in 1975, who considered relativistic
2D axisymmetric problems, and Dorfi \cite{article:Dorfi86} in 1986, who considered
fully three-dimensional flow. These approaches
did not use modern shock-capturing techniques; and therefore,
the computed solutions exhibited strong numerical diffusion.
The first direct\footnote{Constrained
transport methods in the tradition of Evans and Hawley
\cite{article:EvHa88} don't {\it directly} use the magnetic
potential, but as is discussed in T\'oth \cite{article:To00},
a discrete magnetic potential is certainly still in the
background in the definitions of the staggered magnetic
and electric field values.} use of the magnetic potential in the context of
modern shock-capturing schemes was due to Londrillo and 
Del Zanna \cite{article:LoZa00}. Subsequent work on
using a magnetic potential with shock-capturing schemes
includes De Sterck \cite{article:De01b}, Londrillo and Del Zanna
\cite{article:LoZa04}, and Rossmanith \cite{article:Ro04b}.

The main focus  of this work is to extend the 2D
constrained transport method introduced by Rossmanith \cite{article:Ro04b}
to three-dimensions. We begin with an outline of our method that
lists all of the key steps necessary to completely advance the solution from
its current state at time $t=t^n$ to its new state at time $t=t^{n+1}=t^n+\Delta t$:
\begin{description}
  \item[\qquad {\bf Step 0.}]  Start with the current state: 
  \quad $\left( \rho^n, \, \rho\u^n, \, {\mathcal E}^n, \, \B^n, \, \Av^n \right)$.
   \item[\qquad {\bf Step 1.}]  Update MHD variables via 	
   	a standard finite volume scheme:
	\[
	\left( \rho^{n}, \, \rho\u^{n}, \, {\mathcal E}^{n}, \, \B^{n} \right)
	\quad \Longrightarrow \quad
	\left( \rho^{n+1}, \, \rho\u^{n+1}, \, {\mathcal E}^{\star}, \, \B^{\star} \right),
	\]
	where the energy and magnetic field values,  ${\mathcal E}^{\star}$ and 
	$\B^{\star}$, are given a $\star$ superscript instead of $n+1$ to emphasize
	that these values will be modified by the constrained transport
	algorithm before the end of the time step.
   \item[\qquad {\bf Step 2.}]  Define the time-averaged velocity:
   		\quad $\u^{n+\half} = \frac{1}{2} \left( \u^n + \u^{n+1} \right)$.
   \item[\qquad {\bf Step 3.}] Using the above calculated velocity, $\u^{n+\half}$,
   	solve the magnetic potential version of the induction equation\footnote{See
	\S\ref{sec:vecpotential} for an explanation of this equation.}:
		\[   \Av_{,t}  + \left( \nabla \times \Av \right) \times 
			\u^{n+\half} = - \nabla \psi. \]
	This updates the vector potential: $\Av^n  \quad
		 \Longrightarrow  \quad \Av^{n+1}$.
   \item[\qquad {\bf Step 4.}] Compute the new magnetic
   	field from the curl of the vector potential ($\B = \nabla \times \Av$):
	\begin{align}
		\left[ B^1 \right]^{n+1}_{ijk} &= 
			\frac{\left[\A^3\right]^{n+1}_{i\,j+1\,k} - \left[\A^3 \right]^{n+1}_{i\,j-1\,k}}{2 \Delta y}
		     -  \frac{\left[ \A^2 \right]^{n+1}_{i\,j\,k+1} - \left[ \A^2 \right]^{n+1}_{i\,j\,k-1}}{2 \Delta z}, \\
		\left[ B^2 \right]^{n+1}_{ijk} &= \label{eqn:B2}
			\frac{\left[ \A^1 \right]^{n+1}_{i\,j\,k+1} - \left[ \A^1 \right]^{n+1}_{i\,j\,k-1}}{2 \Delta z}
		     -  \frac{\left[ \A^3 \right]^{n+1}_{i+1\,j\,k} - \left[ \A^3 \right]^{n+1}_{i-1\,j\,k}}{2 \Delta x}, \\
		\left[ B^3 \right]^{n+1}_{ijk} &= \label{eqn:B3}
			\frac{\left[ \A^2 \right]^{n+1}_{i+1\,j\,k} - \left[ \A^2 \right]^{n+1}_{i-1\,j\,k}}{2 \Delta x}
		     -  \frac{\left[ \A^1 \right]^{n+1}_{i\,j+1\,k} - \left[ \A^1 \right]^{n+1}_{i\,j-1\,k}}{2 \Delta y},
	\end{align}
	where the bracket notation, $\left[ \cdot \right]^n_{ijk}$, is used to clearly distinguish
	between  vector component superscripts  and grid and time point
	superscripts and subscripts.
  \item[\qquad {\bf Step 5.}] Set the new total energy density 
  value based on one of the following options\footnote{In this paper 
  we always choose {\bf Option 1} in order to exactly conserve energy.}:
    \begin{description}
    	\item[\qquad {\bf Option 1:}] Conserve total energy: \quad 
		\[ {\mathcal E}^{n+1}
		= {\mathcal E}^{\star}. \]
	\item[\qquad {\bf Option 2:}] Keep the pressure the same
		before and after the constrained transport step (this sometimes helps
		in preventing the pressure from becoming negative, although
		it sacrifices energy conservation):
		 \quad \[ {\mathcal E}^{n+1}
		= {\mathcal E}^{\star} + \frac{1}{2} \left( \| \B^{n+1} \|^2 - \| \B^{\star} \|^2
			\right). \]
    \end{description}
\end{description}

The above described algorithm guarantees that at each time step, the following
discrete divergence is identically zero:
\begin{equation}
\begin{split}
	\left[ \nabla \cdot \B \right]^{n+1}_{ijk}  &= 
	\frac{\left[B^1 \right]^{n+1}_{i+1 \, j \, k} - \left[B^1 \right]^{n+1}_{i-1 \, j \, k}}{2 \Delta x} + 
	\frac{\left[B^2 \right]^{n+1}_{i \, j+1 \, k} - \left[B^2 \right]^{n+1}_{i \, j-1 \, k}}{2 \Delta y}  \\
	&+ \, \frac{\left[B^3 \right]^{n+1}_{i \, j \, k+1} - \left[B^3 \right]^{n+1}_{i \, j \, k-1}}{2 \Delta z} = 0.
\end{split}
\end{equation}

\section{Vector potential equations and gauge conditions}
\label{sec:vecpotential}
The key step in the constrained transport framework as outlined
in the previous section is {\bf Step 3}, which requires one to solve
an evolution equation for the magnetic vector potential. One
question that immediately arises: what should be chosen
for the gauge condition? In this section we briefly discuss
several gauge conditions and their consequences
on the evolution of the magnetic vector potential.

The starting point for this discussion is the induction equation:
\be
\label{eqn:induction}
\B_{,t} + \nabla \times \left( \B \times \u \right) = 0,
\ee
where, for the purposes of the algorithm outlined in
\S \ref{sec:framework}, we take the velocity, $\u$,
as a given function.
We set $\B = \nabla \times \Av$ and rewrite \eqref{eqn:induction} as
\begin{gather}
\nabla \times \left\{  \Av_{,t}  + \left( \nabla \times \Av \right) \times \u \right\} = 0, \\
\label{eqn:vecpot_general}
\Longrightarrow \quad \Av_{,t} +  \left( \nabla \times \Av \right) \times \u = -\nabla \psi,
\end{gather}
where $\psi$ is an arbitrary scalar function. Different choices of $\psi$
represent different {\it gauge condition} choices.

Before we explore various gauge conditions, however,
it is worth pointing out that the situation in the pure two-dimensional case
(e.g., in the $xy$-plane) is much simpler.
The only component of the magnetic vector potential
that influences the evolution in this case is
$\A^3$ (i.e., the component of the potential that is
perpendicular to the evolution plane); and furthermore,
 all gauge choices lead to the same equation:
 \be
	\A^3_{,t} + u^1 \A^3_{,x} + u^2 \A^3_{,y} = 0,
\ee
where $\A^3$ is uniquely defined up to an additive constant.

\subsection{Coulomb gauge}
An obvious choice for the gauge is to take the vector magnetic potential
to be solenoidal:
\begin{equation}
	\nabla \cdot \Av = 0,
\end{equation}
resulting in the {\it Coulomb gauge}.
We are now able to add $-	\u \, \nabla \cdot \Av$ to 
the left-hand side of \eqref{eqn:vecpot_general}
and obtain an equation for the potential that
is in symmetric hyperbolic form:
\begin{equation}
\begin{split}
\begin{bmatrix}
  \A^1 \\ \A^2 \\ \A^3
\end{bmatrix}_{,t} 
& +
\begin{bmatrix}
  -u^1 & -u^2 & -u^3 \\
  -u^2 & u^1 & 0 \\
  -u^3 & 0 & u^1
\end{bmatrix}
\begin{bmatrix}
  \A^1 \\ \A^2 \\ \A^3
\end{bmatrix}_{,x}
 +
\begin{bmatrix}
   u^2 & -u^1 & 0 \\
  -u^1 & -u^2 & -u^3 \\
  0 & -u^3 & u^2
\end{bmatrix}
\begin{bmatrix}
  \A^1 \\ \A^2 \\ \A^3
\end{bmatrix}_{,y}
 \\ 
 & + \begin{bmatrix}
  u^3 &  0 &  -u^1 \\
  0 & u^3 & -u^2 \\
  -u^1 &  -u^2 &  -u^3
\end{bmatrix}
\begin{bmatrix}
  \A^1 \\ \A^2 \\ \A^3
\end{bmatrix}_{,z} = -
\begin{bmatrix}
	\psi_{,x} \\  \psi_{,y}  \\  \psi_{,z}
\end{bmatrix}.
\end{split}
\end{equation}
The main difficulty with this approach, however, is
that at each time step one must solve a Poisson equation
to determine the Lagrange multiplier $\psi$:
\begin{equation}
   - \nabla^2 \psi = \nabla \cdot \left[ \left( \nabla \times \Av \right) \times \u \right].
\end{equation}
Having to solve an elliptic equation in each time step makes this approach
have the same efficiency problems as the projection method.

\subsection{Lorentz-like gauge}
In the ideal MHD setting, since the speed of light is taken to be infinite, the Lorentz and Coulomb gauges are equivalent.
However, one possibility is to introduce a fictitious wave speed, $\xi$, that is larger than all other wave
speeds in the MHD system. We can then take
\begin{equation}
	\psi_{,t} = \xi^2 \, \nabla \cdot \Av,
\end{equation}
which results in the following evolution equation for $(\Av, \psi)$:
\begin{equation}
\label{eqn:lorentzgauge}
\begin{split}
\begin{bmatrix}
  \A^1 \\ \A^2 \\ \A^3 \\  \psi
\end{bmatrix}_{,t} 
& +
\begin{bmatrix}
  0 & -u^2 & -u^3 & 1\\
  0 & u^1 & 0 & 0\\
  0 & 0 & u^1 & 0 \\
  \xi^2 & 0 & 0 & 0 \\
\end{bmatrix}
\begin{bmatrix}
  \A^1 \\ \A^2 \\ \A^3 \\ \psi
\end{bmatrix}_{,x}
 +
\begin{bmatrix}
   u^2 & 0 & 0 & 0 \\
  -u^1 & 0 & -u^3 & 1 \\
  0 & 0 & u^2 & 0 \\
  0 & \xi^2 & 0 & 0 
\end{bmatrix}
\begin{bmatrix}
  \A^1 \\ \A^2 \\ \A^3 \\ \psi
\end{bmatrix}_{,y}
 \\ 
 & + \begin{bmatrix}
  u^3 &  0 &  0 & 0 \\
  0 & u^3 & 0 & 0 \\
  -u^1 &  -u^2 &  0 & 1 \\
  0 & 0 & \xi^2 & 0 
\end{bmatrix}
\begin{bmatrix}
  \A^1 \\ \A^2 \\ \A^3 \\ \psi
\end{bmatrix}_{,z} = 0.
\end{split}
\end{equation}
The flux Jacobian of this system in some direction $\n$ (where 
$\| \n \| = 1$) can be written as
\begin{equation}
N(\n) = \begin{bmatrix}
n^2 u^2 + n^3 u^ 3 & -n^1 u^2 & -n^1 u^3 & n^1 \\
-n^2 u^1 & n^1 u^1 + n^3 u^3 & -n^2 u^ 3 & n^2 \\
-n^3 u^1 & -n^3 u^2 & n^1 u^1+n^2 u^2 & n^3 \\
n^1 \xi^2 & n^2 \xi^2 & n^3 \xi^2 & 0
\end{bmatrix}.
\end{equation}
The eigenvalues of this matrix are
\be
  \lambda = \Bigl\{ -\xi, \, \xi, \, \u \cdot \n, \, \u \cdot \n \Bigr\}.
\ee
If $\xi>| \u \cdot \n |$, the right-eigenvectors
of $N(\n)$ are complete, and thus, the system is 
hyperbolic.

Although this seems like a potentially useful gauge 
choice, we found in practice that numerical solutions
to system \eqref{eqn:lorentzgauge} did not produce
accurate magnetic fields. In particular, 
we observed errors in the location of strong shocks,
which is presumably due to the fact that on the discrete
level
\be
\nabla \times \nabla \psi \ne 0, 
\ee
thus resulting in errors in $\B$.

\subsection{Helicity-inspired gauge}
Another gauge possibility is to directly set the scalar function, $\psi$,
to something useful. One such choice is
\begin{equation}
	\psi = \u \cdot \Av,
\end{equation}
which is inspired by the magnetic helicity: $\B \cdot \Av$.
This yields the system:
\begin{align}
\label{eqn:helicity_gauge}
\Av_{,t} 
& +
u^1
\Av_{,x}
 +
u^2
\Av_{,y}
+ 
u^3
\Av_{,z} = M \Av,
\end{align}
where
\begin{equation}
M := -
\begin{bmatrix}
	{u^1_{,x}} & u^2_{,x} & u^3_{,x} \\
	\vspace{-3mm} \\
	u^1_{,y} & u^2_{,y} & u^3_{,y} \\
	\vspace{-3mm} \\
	u^1_{,z} & u^2_{,z} & u^3_{,z} 
\end{bmatrix}.
\end{equation}
One obvious approach for solving this equation is via
operator splitting, whereby equation \eqref{eqn:helicity_gauge} is
split into three decoupled advection equations:
\begin{align}
\Av_{,t} 
& +
u^1
\Av_{,x}
 +
u^2
\Av_{,y}
+ 
u^3
\Av_{,z} = 0,
\end{align}
and a `linear' ordinary differential equation\footnote{This equation is linear
in the sense that the velocity, $\u$, is taken to be frozen in time
at $t=t^{n+1/2}$.}:
\begin{align}
\label{eqn:helicity_ODE}
\Av_{,t}  = M \Av.
\end{align}
The main difficulty with this approach is that
the matrix $M$ in equation \eqref{eqn:helicity_ODE} could
(and often does) have eigenvalues that have a positive real part;
thereby, causing this system to be inherently unstable.
For this reason, numerical tests using this gauge condition
were generally not successful.

\subsection{Weyl gauge}
The choice that we finally settled on was the
Weyl gauge:
\begin{equation}
	\psi = 0.
\end{equation}
In this approach, the resulting evolution equation
is simply \eqref{eqn:vecpot_general} with a zero
right-hand side.
This gauge is the most
commonly used one in the description of
constrained transport methods (see for
example \cite{article:LoZa00,article:To00}).

As we will describe in detail in the next section, \S \ref{sec:numerics},
the resulting system is only {\it weakly} hyperbolic. This is due to the fact
that there are certain directions in which the matrix of right-eigenvectors of the
flux Jacobian does not have full rank. This degeneracy causes some
numerical difficulties, which we were able to overcome through
the creation of a modified wave propagation scheme\footnote{See also
Fey and Torrilhon \cite{article:FeTo03} for a discussion of numerical
discretizations of the weakly hyperbolic induction equation.}. The
details of this approach are described in the next section.

\section{Numerical methods}
\label{sec:numerics}
The numerical methods developed in this paper
are based on the high-resolution wave propagation scheme of
LeVeque \cite{article:Le97} and its 3D extension 
developed by Langseth and LeVeque \cite{article:LaLe00}.
In \S\ref{subsec:LangLev} we briefly review the 3D
wave propagation approach. In \S\ref{sec:method_potential} we show in detail
how this approach can be modified to solve the non-conservative
and only weakly hyperbolic magnetic vector potential equation
\eqref{eqn:system_compact}--\eqref{eqn:N1N2N3}. 
In \S\ref{subsec:mp_limiters} we develop a limiting strategy that is applied
during the magnetic potential update, but designed to control
unphysical oscillations in the magnetic field. Finally, 
in \S\ref{subsec:mp_25d} we briefly mention how our approach simplifies
in the special case of 2.5-dimensional problems.

\subsection{The wave propagation scheme of Langseth and LeVeque
  \cite{article:LaLe00}}
 \label{subsec:LangLev}
In {\bf Step 1} of the constrained transport method described in \S
\ref{sec:framework} we apply a numerical method for the three-dimensional MHD equations.
Here we use a version of the three-dimensional wave propagation
algorithm of Langseth and LeVeque \cite{article:LaLe00}, see also
\cite{book:Le02}, 
which is based on a decomposition of flux differences at grid cell interfaces as
outlined in \cite{article:BaRoLe2}. This is a multidimensional
high-resolution finite volume method that is second order accurate
for smooth solutions and that leads to an accurate capturing of shock
waves.  

To outline the main steps of this algorithm, we consider a
three-dimensional hyperbolic system of the general form 
\begin{equation}
q_{,t} + f(q)_{,x} + g(q)_{,y} + h(q)_{,z} = 0,
\end{equation}
with $q : \mathbb{R}^+ \times \mathbb{R}^3  \rightarrow \mathbb{R}^m$ and
$f, g, h : \mathbb{R}^m \rightarrow \mathbb{R}^m$.
In quasilinear form the system reads
\begin{equation}
q_{, t} + \Am(q) \, q_{, x} + \Bm(q) \, q_{, y} + \Cm(q) \, q_{, z} = {0},
\end{equation}
where $\Am$, $\Bm$, and $\Cm$ are the flux-Jacobians in
each of the three coordinate directions:
\begin{equation}
\Am(q) := f(q)_{,q}, \quad
\Bm(q) := g(q)_{,q}, \quad \text{and} \quad
\Cm(q) := h(q)_{,q}.
\end{equation}
We construct a 3D Cartesian mesh with grid spacing $\Delta x$,
$\Delta y$, and $\Delta z$ in each of the coordinate directions.
We represent the solution at each discrete time level $t=t^n$
as piecewise constant, such that
in each grid cell $(i,j,k)$ the solution is given by $Q_{ijk}^n$.
$Q_{ijk}^n$ denotes the cell average of the conserved quantity
$q(t^n,x,y,z)$ in the grid cell centered at $(x_{i}, y_{j}, z_{k})$:
\be
	Q^{n}_{ijk} \approx 
	\frac{1}{\Delta x \Delta y \Delta z}
	\int_{x_{i}-\frac{\Delta x}{2}}^{x_{i}+\frac{\Delta x}{2}}
	\int_{y_{j}-\frac{\Delta y}{2}}^{y_{j}+\frac{\Delta y}{2}}
	\int_{z_{k}-\frac{\Delta z}{2}}^{z_{k}+\frac{\Delta z}{2}}
	q\left(t^n, \xi, \eta, \zeta \right) \, d\xi \, d\eta \, d\zeta.
\ee
 The numerical method can be
written in the form
\begin{equation}\label{eqn:update3d}
\begin{split}
Q_{ijk}^{n+1} = Q_{ijk}^n & - \frac{\Delta t}{\Delta x} \left(
  {\cal A} ^+ \Delta Q_{i-\frac{1}{2} \, j \,k} + {\cal A}^- \Delta
  Q_{i+\frac{1}{2} \, j \, k} \right) \\
& - \frac{\Delta t}{\Delta y} \left( {\cal B}^+ \Delta
  Q_{i \, j-\frac{1}{2} \, k} + {\cal B}^- \Delta Q_{i \, j+\frac{1}{2} \,k}
\right) \\
& - \frac{\Delta t}{\Delta z} \left( {\cal C}^+ \Delta
  Q_{i \, j \, k-\frac{1}{2}} + {\cal C}^- \Delta Q_{i \, j \, k+\frac{1}{2}}
\right)\\
& - \frac{\Delta t}{\Delta x} \left( \tilde F_{i+\frac{1}{2} \, j \, k} -
  \tilde F_{i-\frac{1}{2} \,j \,k} \right) - \frac{\Delta t}{\Delta y}
\left( \tilde G_{i \, j+\frac{1}{2} \, k} -
  \tilde G_{i \, j-\frac{1}{2} \, k} \right) \\
&- \frac{\Delta t}{\Delta z} \left( \tilde H_{i \, j \, k+\frac{1}{2}} -
  \tilde H_{i \, j \, k-\frac{1}{2}} \right). 
\end{split}
\end{equation}
The first three lines in (\ref{eqn:update3d}) describe a first order
accurate update and the last two lines represent higher order
correction terms. The wave propagation method is a so-called {\it truly
multidimensional scheme} in the sense that no dimensional splitting is
used to approximate the mixed derivative terms that are
required in  a second order accurate update.

At each grid cell interface in the $x$-direction, we decompose the
flux differences
\begin{equation}
\Delta F_{i-\frac{1}{2} \, j \, k} = f (Q_{i \, j \, k}^n) - f(Q^n_{i-1 \, j \, k}),
\end{equation} 
into waves 
\begin{equation}
{\cal Z}_{i-\frac{1}{2} \, j \, k}^p = \left[ \l_{i-\frac{1}{2} \, j \, k}^p \cdot
\Delta F_{i-\frac{1}{2} \, j \, k} \right] \, \r_{i-\frac{1}{2} \, j \, k}^p, 
\end{equation}
which are moving with speeds $s_{i-\frac{1}{2} \, j \, k}^p$, $p=1, \ldots, M_w$.
For this decomposition we use the left and right eigenvectors
proposed by Roe and Balsara \cite{article:RoBa96} (see also Powell \cite{article:Po99}
for a discussion of these eigenvectors) for a 
linearized flux Jacobian matrix of the MHD equations.
The fluctuations used in the first order update have the form:
\begin{align}
{\cal A}^+ \Delta Q_{i-\frac{1}{2} \, j \, k} & = \sum_{p:
  s_{i-\frac{1}{2} \, j \, k}^p > 0} {\cal Z}_{i-\frac{1}{2} \, j \, k}^p +
\frac{1}{2} \sum_{p: s_{i-\frac{1}{2} \, j \, k}^p = 0} {\cal
  Z}_{i-\frac{1}{2} \, j \, k}^p, \\
{\cal A}^- \Delta Q_{i-\frac{1}{2} \, j \, k} & = \sum_{p:
  s_{i-\frac{1}{2} \, j \, k}^p < 0} {\cal Z}_{i-\frac{1}{2} \, j \, k}^p +
\frac{1}{2} \sum_{p: s_{i-\frac{1}{2} \, j \, k}^p = 0} {\cal
  Z}_{i-\frac{1}{2} \, j \, k}^p.
\end{align}
Analogously, we compute fluctuations in the $y-$ and $z-$directions.
An advantage of the decomposition of the flux differences compared to
the standard wave propagation method from \cite{article:LaLe00,article:Le97},
 which is based on an eigenvector
decomposition of the jumps of the conserved quantities at each grid cell
interface, is that the local linearization of the flux Jacobian matrix
does not require Roe averages, which are quite difficult to compute for
the MHD equations \cite{article:CaCa97}.
Here we use simple arithmetic averaging instead.
See Wesenberg \cite{article:Wesenberg02} for a discussion
of how the use of arithmetic averages is equivalent to an MHD Riemann
solver based on van Leer flux-vector splitting \cite{article:vanLeer82}.

The waves and speeds are also used to compute high-resolution
correction terms. For the waves in the $x$-direction, the correction
fluxes have
the form
\begin{equation}
\label{eqn:limiter1}
\tilde F_{i-\frac{1}{2} \, j \, k} = \frac{1}{2} \sum_{p=1}^{M_w}
\mbox{sign}\left(s_{i-\frac{1}{2} \, j \, k}^p\right) \left( 1 - \frac{\Delta t}{\Delta x}
  \left| s_{i-\frac{1}{2} \, j \, k}^p \right| \right) {\cal Z}_{i-\frac{1}{2} \, j \, k}^p \, \phi
   \left( \theta^{p}_{i-\half \, j \, k} \right),
\end{equation}
where 
\be
\label{eqn:limiter2}
	\theta^p_{i-\half \, j \, k} := \frac{{\cal Z}_{i-\frac{1}{2} \, j \, k}^p \cdot {\cal Z}_{I-\frac{1}{2} \, j \, k}^p}{{\cal Z}_{i-\frac{1}{2} \, j \, k}^p \cdot {\cal Z}_{i-\frac{1}{2} \, j \, k}^p},
	\quad \text{where} \quad I = \begin{cases}
			i-1 & \quad \text{if} \quad s^p_{i-\half \, j \, k} > 0, \\
			i+1 & \quad \text{if} \quad  s^p_{i-\half \, j \, k} < 0, 
		\end{cases}
\ee
is a measure of the location smoothness in the $p$-th wave
 family and $\phi(\theta)$ is a total variation diminishing (TVD) limiter function
 \cite{article:Ha83, article:Vl74,article:Sw84}.

The wave propagation method is a 1-step, unsplit, second-order
accurate method in both space and time (for smooth solutions);
and therefore, is based on a Taylor series expansions in time:
\be
\begin{split}
	 q \left(t+\Delta t, {\bf x} \right) = q  +	
		\Delta t  q_{,t} + \frac{1}{2} \Delta t^2  q_{,tt} + {\mathcal O} \left( \Delta t^3 \right) 
		\hspace{18mm} \\
	= q  - \Delta t \Bigl[ f_{,x} + g_{,y} + h_{,z} \Bigr] 
		 + \frac{1}{2} \Delta t^2 
		\Bigl\{ \left( \Am  f_{,x} \right)_{,x} + \left( \Bm   g_{,y} \right)_{,y} 
				+  \left( \Cm   h_{,z} \right)_{,z} \hspace{4mm} \\
	 +  \underset{\text{Transverse terms}}{\underbrace{\left( \Am   g_{,y} \right)_{,x} + \left( \Am    h_{,z} \right)_{,x} +
	\left( \Bm    f_{,x} \right)_{,y} + \left( \Bm    h_{,z} \right)_{,y} +
	\left( \Cm   f_{,x} \right)_{,z} + \left( \Cm   g_{,y} \right)_{,z}  }}  \Bigr\} 
	 + {\mathcal O} \left( \Delta t^3 \right),
\end{split}
\ee
where time derivatives have been replaced by spatial derivations through
the conservation law.
The wave propagation method as described so far, approximates each of
the terms in the above Taylor series to second order accuracy,
except those terms marked as {\it transverse terms}. The transverse terms
are approximated via additional Riemann solvers known as {\it transverse
Riemann solvers}. Additionally, Langseth and LeVeque \cite{article:LaLe00}
found that in order to achieve optimal stability bounds, additional so-called
{\it double transverse Riemann solvers} are required to approximate
certain mixed third-derivative terms. We omit a full discussion of the
transverse and double transverse terms in this work, and instead,
refer the reader to the paper of Langseth and LeVeque \cite{article:LaLe00}.

\subsection{The evolution of the magnetic  potential}
\label{sec:method_potential}
We now describe a numerical method for the evolution equation of the
magnetic potential using the Weyl gauge, i.e., for the discretization
of (\ref{eqn:vecpot_general}) with zero right hand side.  The
equation can be written in the form
\begin{equation}\label{eqn:system_compact}
\Av_{, t} + N_1(\u) \, \Av_{,x} + N_2(\u) \, \Av_{, y} + N_3(\u) \, \Av_{, z} = 0,
\end{equation}
with 
\begin{equation}\label{eqn:N1N2N3}
N_1 = \begin{bmatrix}
0 & -u^2 & -u^3\\
0 & u^1 & 0\\
0 & 0 & u^1\end{bmatrix},  \,
N_2 = \begin{bmatrix}
u^2 & 0 & 0\\
-u^1 & 0  & -u^3\\
0 & 0 & u^2\end{bmatrix}, \,
N_3 = \begin{bmatrix}
u^3 & 0  & 0\\
0 & u^3 & 0\\
-u^1 & -u^2 & 0 \end{bmatrix}.
\end{equation}
\subsubsection{Weak hyperbolicity}
We recall the following definitions regarding hyperbolicity.

\bigskip

\noindent
{\bf Definition.} ({\bf Strict hyperbolicity})
{\it The quasilinear system,
\begin{equation}
\label{eqn:quasilinear3d}
\q_{, t} + \Am(\q) \, \q_{, x} + \Bm(\q) \,  \q_{, y} + \Cm(\q) \,  \q_{, z} = {0},
\end{equation}
is  strictly hyperbolic if the matrix
\begin{equation}
\label{eqn:fluxJ}
	M({\bf n},q) := n^1 \, \Am(q) + n^2 \,  \Bm(q) + n^3  \, \Cm(q)
\end{equation}
is diagonalizable with distinct real eigenvalues for all $\| {\bf n} \| = 1$.
In other words, the $p^{\text{th}}$ eigenvalue of 
$M({\bf n},q)$ is real and has  geometric and algebraic multiplicities of
exactly one.}

\bigskip

\noindent
{\bf Definition.} ({\bf Non-strict hyperbolicity})
{\it The quasilinear system \eqref{eqn:quasilinear3d} is 
non-strictly hyperbolic if the matrix \eqref{eqn:fluxJ}
is diagonalizable with real but not necessarily distinct
eigenvalues for all $\| {\bf n} \| = 1$.
In other words, the $p^{\text{th}}$ eigenvalue of 
$M({\bf n},q)$ is real and has the same geometric and algebraic
multiplicity, but that multiplicity may be greater than one.}

\bigskip

\noindent
{\bf Definition.} ({\bf Weak hyperbolicity})
{\it The quasilinear system \eqref{eqn:quasilinear3d} is 
weakly hyperbolic if the matrix \eqref{eqn:fluxJ}
has real eigenvalues for all $\| {\bf n} \| = 1$, but is 
not necessarily diagonalizable.
In other words, the $p^{\text{th}}$ eigenvalue of 
$M({\bf n},q)$ is real but may have an algebraic multiplicity larger
than its geometric multiplicity.}

\bigskip

The flux Jacobian
of system (\ref{eqn:system_compact}) in an arbitrary 
direction $\n \in S^2$ is
\begin{equation}
n^1 N_1 + n^2 N_2 + n^3 N_3 = \begin{bmatrix}
n^2 u^2 + n^3 u^3 & - n^1 u^2 & - n^1 u^3 \\
-n^2 u^1 & n^1 u^1 + n^3 u^3 & - n^2 u^3 \\
-n^3 u^1 & -n^3 u^2 & n^1 u^1 + n^2 u^2
\end{bmatrix}.
\end{equation}
The eigenvalues of this system are 
\be
\lambda = \bigl\{0, \n \cdot \u, \n \cdot \u\bigr\};
\ee
and therefore, we always
have real eigenvalues. The eigenvectors can be written in
the following form:
\begin{equation}
R = \Biggl[ \r^{(1)} \, \Biggl| \, \r^{(2)} \, \Biggl| \, \r^{(3)} \Biggr] = 
\begin{bmatrix}
n^1 & \quad  n^2 u^3 - n^3 u^2 & \quad u^1 \left( \u \cdot \n \right) - n^1 \| \u \|^2 \\
n^2 & \quad n^3 u^1 - n^1 u^3 & \quad u^2 \left( \u \cdot \n \right) - n^2 \| \u \|^2 \\
n^3 & \quad n^1 u^2 - n^2 u^1 & \quad u^3 \left( \u \cdot \n \right) - n^3 \| \u \|^2
\end{bmatrix}.
\end{equation}
Assuming that $\| \u \| \ne 0$ and $\| \n \| = 1$, the determinant of
matrix $R$ can be written as
\begin{equation}
\text{det}(R) = -  \| \u \|^3 \, \cos(\alpha) \, \sin^2(\alpha),
\end{equation}
where $\alpha$ is the angle between the vectors $\n$ and $\u$.
The difficulty is that for any non-zero velocity vector, ${\bf u}$, 
one can always find four directions, $\alpha=0$, $\pi/2$, $\pi$, and $3\pi/2$,
such that $\text{det}(R) = 0$.
In other words, for every $\| {\bf u} \| \ne 0$ there exists four degenerate directions in which the eigenvectors are incomplete.
Therefore, system \eqref{eqn:system_compact}--\eqref{eqn:N1N2N3}
 is only {\em weakly hyperbolic}.

\subsubsection{An example: the difficulty with weakly hyperbolic systems}
Non-strict hyperbolicity is common in many 
standard equations (e.g., Euler equations
of gas dynamics); and, although it causes some
difficulties in proving long-time existence results, 
it generally does not cause problems for numerical
discretization of such equations.
Weak hyperbolicity, however, is a different story.
We illustrate this point
with the following simple example.
Let
\begin{equation}
  \begin{bmatrix}
   u \\ v 
   \end{bmatrix}_{,t} + 
   \begin{bmatrix}
   -\varepsilon & 1 \\
   0 & \varepsilon
   \end{bmatrix}
   \begin{bmatrix}
   u \\ v 
   \end{bmatrix}_{,x} = 0,
\end{equation}
where $\varepsilon \in \reals$ is a constant. 
The eigen-decomposition of the flux Jacobian can be written
as follows:
\begin{equation}
\begin{bmatrix}
   -\varepsilon & 1 \\
   0 & \varepsilon
   \end{bmatrix} = 
   R \Lambda R^{-1} =
   \begin{bmatrix}
   1 & 1 \\
   0 & 2\varepsilon
   \end{bmatrix}
   \cdot
   \begin{bmatrix}
   -\varepsilon & 0 \\
   0 & \varepsilon
   \end{bmatrix}
   \cdot
   \frac{1}{2\varepsilon}
\begin{bmatrix}
   2\varepsilon & -1 \\
   0 & 1
   \end{bmatrix}.
\end{equation}
Since the eigenvalues are always real, this system is hyperbolic.
For all $\varepsilon \ne 0$, the system is strongly hyperbolic, and for
$\varepsilon = 0$ the system is only weakly hyperbolic.
Since we have the exact eigen-decomposition, we can write down
the exact solution for the Cauchy problem for all $\varepsilon$:
\begin{equation}
 \begin{bmatrix}
   u \\ v 
   \end{bmatrix} = 
   \begin{bmatrix}
   u_0(x + \varepsilon t) - \frac{1}{2\varepsilon} \bigl\{ v_0(x+\varepsilon t) -
   	v_0(x-\varepsilon t) \bigr\} \\ v_0(x-\varepsilon t)
   \end{bmatrix}.
\end{equation}
In the weakly hyperbolic limit we obtain:
\begin{equation}
 \lim_{\varepsilon \rightarrow 0}
 \begin{bmatrix}
   u \\ v 
   \end{bmatrix} = 
   \begin{bmatrix}
   u_0(x) -  t  \, v'_0(x) \\
   v_0(x)
   \end{bmatrix}.
\end{equation}
The change in dynamics between the strongly and weakly hyperbolic
regimes is quite dramatic: in the strongly hyperbolic case the total
energy is conserved, while in the weakly hyperbolic case the
amplitude of the solution grows linearly in time.

\subsubsection{An operator split approach}
In the case of the magnetic vector potential equation \eqref{eqn:system_compact}--\eqref{eqn:N1N2N3},
the weak hyperbolicity is only an artifact of how we are evolving the magnetic
potential. One should remember that the underlying equations (i.e., the ideal
MHD system) are in fact hyperbolic. One way to understand all of this is
to remember that the true velocity field, $\u$, depends on the magnetic
potential in a nonlinear way. That is, even though there might be instantaneous
{\it growth} in the vector potential due to the weak hyperbolicity in the potential
equation, this growth will immediately change the velocity field in such a way 
that the overall system remains hyperbolic.

Numerically, however, we must construct a modified finite volume method
that can handle this {\it short-lived} weak hyperbolicity. Through some experimentation with
various ideas, we found that simple operator splitting techniques lead to robust
numerical methods. In particular, we found that there are two obvious ways to split system
(\ref{eqn:vecpot_general}) into sub-problems. 

The first possibility is based on a decomposition of the form
\begin{gather}
\label{eqn:splitting1a}
\begin{split}
\mbox{{\bf Sub-problem 1:} \hspace*{0.5cm}} & \A^1_{, t} + u^2 \A^1_{, y} + u^3
\A^1_{, z} = 0,\\
& \A^2_{, t} - u^1 \A^1_{, y} = 0,\\
& \A^3_{, t} - u^1 \A^1_{, z}  = 0,\\[0.5cm]   
\end{split} \\
\label{eqn:splitting1b}
\begin{split}
\mbox{{\bf Sub-problem 2:} \hspace*{0.5cm}} & \A^1_{, t} - u^2 \A^2_{, x} = 0,\\
& \A^2_{, t} + u^1 \A^2_{, x} + u^3 \A^2_{, z} = 0,\\
& \A^3_{, t} - u^2 \A^2_{, z}  = 0, \\[0.5cm]
\end{split}\\
\label{eqn:splitting1c}
\begin{split}
\mbox{{\bf Sub-problem 3:} \hspace*{0.5cm}} & \A^1_{, t} - u^3 \A^3_{, x} = 0,\\
& \A^2_{, t} - u^3 \A^3_{, y} = 0,\\
& \A^3_{, t} + u^1 \A^3_{, y} + u^2 \A^3_{, x} = 0.   
\end{split}
\end{gather}
First, consider Sub-problem 1. The first equation can be approximated
using the two-dimensional method described in \cite[Section
5.3]{article:Ro04b}, which is a modification of LeVeque's wave
propagation method  that leads to a second order
accurate approximation for smooth solutions, both in the
potential and the magnetic field, and a non-oscillatory
solution simultaneously for both the magnetic potential and magnetic field.
Once the quantity $\A^1$ is updated, central finite difference
approximations can be used to update $\A^2$ and $\A^3$.  
Sub-problems 2 and 3 can be approached analogously and
Strang operator splitting can be used to approximate (\ref{eqn:system_compact}).
We have tested this splitting and obtained very satisfying results. 

A second approach to split (\ref{eqn:vecpot_general}) is based on 
dimensional splitting, i.e., we consecutively solve the problems
\begin{gather}
\label{eqn:splitting2a}
\begin{split}
\mbox{{\bf Sub-problem 1:} \hspace*{0.5cm}} & \A^1_{, t} - u^2 \A^2_{, x} -
u^3 \A^3_{, x} = 0,\\
& \A^2_{, t} + u^1 \A^2_{, x} = 0,\\
& \A^3_{, t} + u^1 \A^3_{, x} = 0,\\[0.5cm]
\end{split} \\
\label{eqn:splitting2b}
\begin{split}
\mbox{{\bf Sub-problem 2:} \hspace*{0.5cm}} & \A^1_{, t} + u^2 \A^1_{, y} = 0,\\
& \A^2_{, t} - u^1 \A^1_{, y} - u^3 \A^3_{, y} = 0,\\
& \A^3_{, t} + u^2 \A^3_{, y} = 0,\\[0.5cm]
\end{split} \\
\label{eqn:splitting2c}
\begin{split}
\mbox{{\bf Sub-problem 3:} \hspace*{0.5cm}} & \A^1_{, t} + u^3 \A^1_{, z} =
0,\\
& \A^2_{, t} + u^3 \A^2_{, z} = 0,\\
& \A^3_{, t} - u^1 \A^1_{, z} - u^2 \A^2_{, z} = 0.
\end{split}
\end{gather}
We will now present a numerical method that  is based on this
second splitting. Our method should be easy to adapt by other users; and,  in
particular, also by those who do not base their numerical method on the wave 
propagation algorithm. 

We denote the numerical solution operator for Sub-problem 1, 2 and 3 by
$L_1$, $L_2$ and $L_3$, respectively. 
Then a Strang-type operator splitting method can be written in the form
\be
\Av^{n+1} = L_1^{\Delta t/2} L_2^{\Delta t/2} L_3^{\Delta t} L_2^{\Delta t/2} L_1^{\Delta t/2} \Av^n.
\ee
One time step from $t^n$ to $t^{n+1}$ consists in consecutively solving
Sub-problem 1 over a half time step, Sub-problem 2 over a
half time step, Sub-problem 3 over a full time step, Sub-problem
2 over another half time step and finally solving Sub-problem 1 over another
half time step.   

\subsubsection{Discretization of Sub-problem 1}
\label{sec:subprob1}
We now present the details for the approximation of the first
Sub-problem (\ref{eqn:splitting2a}) 
for a time step from $\tau^n$ to $\tau^{n+1}$. The evolution of $\A^2$ and $\A^3$ is described by two
decoupled scalar transport equations, which have the form
\begin{align}
\label{eqn:sub1_A2}
 \A^2_{,t} + u^1(\tau^{n+\frac{1}{2}},{\bf x}) \,  \A^2_{,x} &=
0,\\
\label{eqn:sub1_A3}
\A^3_{,t} + u^1(\tau^{n+\frac{1}{2}},{\bf x}) \, \A^3_{,x} &=
0.
\end{align}
We assume that cell average values of $\A^{2}$ and $\A^{3}$ are given at
time $\tau^n$ and we wish to approximate $\A^{2}$ and $\A^{3}$ at time $\tau^{n+1}$.
In our application, cell centered values of the advection speed $u^1$ are 
given at the intermediate time level (see {\bf Step 2} of the
algorithm from \S \ref{sec:framework}). 

The update for all $i, j, k$ and $\l \in \{2,3\}$ has the form
\begin{equation}
\begin{split}
\bigl[\A^\l \bigr]_{ijk}^{n+1} = \bigl[ \A^\l \bigr]_{ijk}^n & - \frac{\Delta \tau}{\Delta x} \left[ {\cal
    A}^- \Delta \A^\l_{i+\frac{1}{2} \, j \, k} + {\cal A}^+ \Delta
  \A^\l_{i-\frac{1}{2} \, j \, k} \right] \\
& - \frac{\Delta \tau}{\Delta x} \left[
  \tilde{F}_{i+\frac{1}{2} \, j \, k}^{\ell-}  - \tilde{F}_{i-\frac{1}{2} \, j \, k}^{\ell+} \right],
\end{split}
\end{equation}
with
\begin{align}
{\cal W}^\l_{i-\frac{1}{2} \, j \, k} & := \bigl[\A^\l \bigr]_{ijk}^n - \bigl[\A^\l \bigr]_{i-1\, j \, k}^n,\\
{\cal A}^- \Delta \A^\l_{i-\frac{1}{2} \, j \, k} & := \min(u_{i-1 \, j \, k}^1,0) \, {\cal
  W}^\l_{i-\frac{1}{2} \, j \, k}, \\
{\cal A}^+ \Delta \A^\l_{i-\frac{1}{2} \, j \, k} & := \max(u_{i \, j \, k}^1,0) \, {\cal
  W}^\l_{i-\frac{1}{2} \, j \, k},
\end{align}
and 
\begin{align}
\tilde{F}_{i+\frac{1}{2} \, j \, k}^{\ell-}  & := \frac{1}{2} |u_{ijk}^1| \left( 1
  - \frac{\Delta \tau}{\Delta x}
  \frac{u_{i+\frac{1}{2}\,j\,k}^1}{|u_{ijk}^1|} \, u_{ijk}^1 \right)
{\cal W}^\l_{i+\frac{1}{2}\,j\,k} \, \phi\left(\theta_{ijk}^{\ell} \right),\\
\tilde{F}_{i-\frac{1}{2}\,j\,k}^{\ell+}  & := \frac{1}{2} |u_{ijk}^1| \left( 1 -
  \frac{\Delta \tau}{\Delta x} \frac{u_{i-\frac{1}{2}\,j\,k}^1}{|u_{ijk}^1|}
   \, u_{ijk}^1 \right) {\cal W}^\l_{i-\frac{1}{2}\,j\,k} \,
  \phi\left(\theta_{ijk}^{\ell} \right).
\end{align}
Here $u_{ijk}^1$ is a cell centered value and
$u_{i-\frac{1}{2}jk}^1 = \frac{1}{2}
(u_{i-1jk}^1+u_{ijk}^1)$ is a face centered value.
The key innovation in this construction as developed
in \cite{article:Ro04b} is the use of a cell-based limiter,
rather than the standard edge-based limiter
as represented in equations \eqref{eqn:limiter1}.
Furthermore, in this construction the smoothness
indicator, $\theta_{ijk}^{\ell}$, is based on {\it wave-differences}
rather than waves as in \eqref{eqn:limiter2}:
\be
  \theta_{ijk}^{\ell} = \frac{\Delta {\cal W}^{\ell}_{Ijk}}{\Delta
  {\cal W}^\l_{ijk}},
\ee
where $\Delta {\cal W}^{\ell}_{ijk} = {\cal W}^{\ell}_{i+\frac{1}{2}jk} 
- {\cal W}^{\ell}_{i-\frac{1}{2}jk}$, and the index $I$ is again chosen from
the upwind direction. The reason for the use of wave-differences
instead of waves is that ultimately, the physical quantity of importance
is the magnetic field and not the magnetic potential.
As was shown in \cite{article:Ro04b}, standard wave limiters
do not control oscillations in derivative quantities (i.e., the magnetic
field), while the wave-difference limiter approach does.

Once we have updated $\A^2$ and $\A^3$, we update $\A^1$ 
by discretizing the equation
\begin{equation}
\label{eqn:sub1_A1}
\A^1_{, t} - u^2 \A^2_{, x} - u^3 \A^3_{, x} = 0.
\end{equation}
For this we use a 
finite difference method of the form
\begin{equation}
\label{eqn:weakhyp_solve}
\begin{split}
\left[\A^1\right]_{ijk}^{n+1} = \left[\A^1\right]_{ijk}^n + \frac{\Delta \tau}{2} \biggl\{ &
  \left[u^2\right]_{ijk}^{n+\half} \left(D_x \left[\A^2\right]_{ijk}^n +  D_x
  \left[\A^2\right]_{ijk}^{n+1} \right) \\
+ &  \left[u^3\right]_{ijk}^{n+\half} \left( D_x \left[\A^3\right]_{ijk}^n + D_x
  \left[\A^3\right]_{ijk}^{n+1} \right) \biggr\},
\end{split}
\end{equation}
where $D_x$ is the standard centered finite difference operator in the
$x$-direction, i.e.,
\be
D_x \left[\A^2\right]_{ijk}^n := \frac{
  \left[\A^2\right]_{i+1jk}^n-\left[\A^2\right]_{i-1jk}^n}{2 \Delta x}.
\ee
We will refer to the numerical solution of  \eqref{eqn:sub1_A2} and \eqref{eqn:sub1_A3}
as {\it hyperbolic solves}  and the numerical solution of 
equation \eqref{eqn:sub1_A1} as the {\it weakly hyperbolic solve}.

\subsection{Additional limiting in the weakly hyperbolic solve}
\label{subsec:mp_limiters}
For the
spatially two-dimensional case, the limiting of wave-differences 
eliminates unphysical oscillations in the magnetic field.
In the three-dimensional case, the
equation for the magnetic potential is more challenging 
and only weakly hyperbolic. Each sub-problem of the
numerical method is a combination of one-dimensional hyperbolic solves
(i.e., equations \eqref{eqn:sub1_A2} and \eqref{eqn:sub1_A3}) and
a part that we refer to as the weakly hyperbolic solve (i.e., equation
\eqref{eqn:sub1_A1}). 
The wave-difference limiting strategy as outlined in \S\ref{sec:subprob1}
is generally insufficient in controlling unphysical oscillations in
the magnetic field variables. In particular, the step that
produces unphysical oscillations is the weakly hyperbolic solve
as described in equation \eqref{eqn:weakhyp_solve}. Since this part of the update
is quite different from standard upwind numerical methods for hyperbolic
equations, there is no obvious place to introduce wave-difference limiters.
Instead, we describe in this section an approach based on adding
artificial diffusion to this part of the update in order to remove these unphysical oscillations.

We introduce a diffusive limiter
that is inspired by the artificial viscosity method
that is often used in other numerical schemes such as
the discontinuous Galerkin approach (e.g., see \cite{article:PerPer06}).
Instead of (\ref{eqn:sub1_A1}), we discretize the problem with
an added dissipative term
\begin{equation}
\A^1_{, t} - u^2 A^2_{,x} - u^3 A^3_{,x} 
= \varepsilon^1 A^1_{,xx},
\end{equation}  
where the parameter $\varepsilon^1$ depends on the solution structure
and controls the amount of artificial diffusion. We choose 
\begin{equation}
\varepsilon^1 = 2 \nu^1 \alpha^1 \frac{\Delta x^2}{\Delta t},
\end{equation}
where $\nu^1$ is a positive constant that will be discussed below,
$\alpha^1$
is a smoothness indicator that is close to zero in smooth regions and
close to 0.5 near discontinuities.
Note that we distinguish the size of the total time step:
\be
	\Delta t = t^{n+1} - t^n,
\ee
from the time step of a particular substep of the Strang splitting:
\be
	\Delta \tau = \tau^{n+1} - \tau^n.
\ee
We add this additional diffusion term to each weakly hyperbolic
equation in our dimensionally split algorithm in the following way:
\begin{gather}
\label{eqn:splitting_viscous_a}
\begin{split}
\mbox{{\bf Sub-problem 1:} \hspace*{0.5cm}} & \A^1_{, t} - u^2 \A^2_{, x} -
u^3 \A^3_{, x} = \varepsilon^1 A^1_{,xx},\\
& \A^2_{, t} + u^1 \A^2_{, x} = 0,\\
& \A^3_{, t} + u^1 \A^3_{, x} = 0,\\[0.5cm]
\end{split} \\
\label{eqn:splitting_viscous_b}
\begin{split}
\mbox{{\bf Sub-problem 2:} \hspace*{0.5cm}} & \A^1_{, t} + u^2 \A^1_{, y} = 0,\\
& \A^2_{, t} - u^1 \A^1_{, y} - u^3 \A^3_{, y} = \varepsilon^2 A^2_{,yy},\\
& \A^3_{, t} + u^2 \A^3_{, y} = 0,\\[0.5cm]
\end{split} \\
\label{eqn:splitting_viscous_c}
\begin{split}
\mbox{{\bf Sub-problem 3:} \hspace*{0.5cm}} & \A^1_{, t} + u^3 \A^1_{, z} =
0,\\
& \A^2_{, t} + u^3 \A^2_{, z} = 0,\\
& \A^3_{, t} - u^1 \A^1_{, z} - u^2 \A^2_{, z} = \varepsilon^3 A^3_{,zz}.
\end{split}
\end{gather}
Note that we diffuse only in the direction of the dimensionally split
solve.

Consider again only the discretization
used in Sub-problem (\ref{eqn:splitting_viscous_a}).
The numerical update in the weakly hyperbolic solve can
be written as 
\begin{equation}
\begin{split}
\left[ \A^1 \right]_{ijk}^{n+1} 
& = \left[ \A^1 \right]_{ijk}^{\star} + 2 \frac{\Delta \tau}{\Delta t} \nu \alpha \left(
\left[ \A^1 \right]_{i-1jk}^n - 2 \left[ \A^1 \right]_{ijk}^n + \left[ \A^1 \right]_{i+1jk}^n \right),
\end{split}
\end{equation}
where $\left[ A^1 \right]_{ijk}^{\star}$ is given from update \eqref{eqn:weakhyp_solve}.
Stability requires that
\be
0 \le 2 \frac{\Delta \tau}{\Delta t} \nu \alpha \le \frac{1}{2},
\ee
which can be achieved if we assume that $\Delta \tau \le \Delta t$ and
we take
\be
   0 \le \alpha \le \frac{1}{2} \qquad \text{and} \qquad
   0 \le \nu \le \frac{1}{2}.
\ee
In our proposed method, $\alpha$ is a smoothness indicator and
 $\nu$ is a user-defined parameter (typically $\nu \ll 1$) that can
 be tuned to control the amount of numerical dissipation across
discontinuities in the magnetic field.

The smoothness indicator $\alpha$ is computed according to the formulas
\begin{equation}
\alpha = \max \left( \left| \frac{a_{\text {l}}}{a_{\text {l}}+a_{\text {r}}} - \half \right|, \left|
    \frac{a_{\text {r}}}{a_{\text {l}}+a_{\text {r}}} - \half \right| \right),
\end{equation}
with
\begin{equation}
a_{\text {l}} = \left\{ {\epsilon} + \left( \A^{1}_{ijk} - \A^{1}_{i-1jk}
  \right)^2 \right\}^{-2}, \
a_{\text {r}} = \left\{ {\epsilon} + \left( \A^{1}_{i+1jk} - \A^{1}_{ijk}
  \right)^2 \right\}^{-2}.
\end{equation}
The parameter $\epsilon$ is introduced only to avoid division theory and
is in practice taken to be ${\epsilon} = 10^{-8}$. The smoothness indicator is designed
to distinguish between the following cases:
\begin{enumerate}
\item If $\A^1$ is smooth for $x \in (x_{i-1}, x_{i+1})$, then
	 \begin{align*}
	 \A^{1}_{ijk} - \A^{1}_{i-1jk} = {\mathcal O}(\Delta x), \qquad \text{and} \qquad
	  \A^{1}_{i+1 \, jk} - \A^{1}_{ijk}  =  {\mathcal O}(\Delta x).
	  \end{align*}
	  In this case one can show that
	  \[
	  	\lim_{\Delta x \rightarrow 0}
		\lim_{\epsilon \rightarrow 0}
		\alpha = \left| \frac{A^1_{,xx}}{A^1_{,x}} \right|
			\Delta x = {\mathcal O}(\Delta x),
	  \]
	  which shows that the overall numerical method still retains
	  ${\mathcal O}(\Delta t^2,\Delta x^2)$ accuracy.
\item If $\A^1$ is non-smooth in $(x_{i-1}, x_{i})$ and
smooth in $(x_{i}, x_{i+1})$, then
	$a_\text{r} \gg a_\text{l}$ and $\alpha \approx \frac{1}{2}$.
\item If $\A^1$ is non-smooth in  $(x_{i}, x_{i+1})$ and smooth
in $(x_{i-1}, x_{i})$, then 
         $a_\text{l} \gg a_\text{r}$ and $\alpha \approx \frac{1}{2}$.
\end{enumerate}
 In all cases $\alpha \le 1/2$, which guarantees that the numerical
 update will be stable up to CFL number one.

\subsection{Discretization of the 2.5-dimensional problem}
\label{subsec:mp_25d}
In order to construct and analyze methods for the three-dimensional
MHD equations, it is instructive to also look at the so-called 2.5-dimensional 
case. In this case we assume that $\u$ and $\B$ are
three-dimensional vectors, but all conserved quantities
$[\rho, \rho \u, {\cal E}, \B]^T$ are functions of only two spatial
variables, say $\x = (x,y)^T$. 

Since $B^3_{, z} = 0$, we can
obtain a divergence free magnetic field by employing a constrained
transport algorithm that only updates $B^1$ and $B^2$, and which treats $B^3$
as just another conserved variable in the base scheme.
In other words, one could simply update the transport
equation for $A^3$,
\begin{equation}\label{eqn:mp3}
\A^3_{, t} + u^1 \A^3_{, x} + u^2 \A^3_{, y} = 0,
\end{equation}
via the method of \cite{article:Ro04b}, then compute
the first two components of the magnetic field from 
discrete analogs of  $B^1 = A^3_{,y}$ and $B^2=-A^3_{,x}$, and
let the base scheme update $B^3$ since it doesn't enter
into the 2.5D divergence constraint:
\[
	B^1_{,x} + B^2_{,y} = 0.
\]

An alternative approach to handling the 2.5D case, one that will allow
us to test the important features of the proposed
3D scheme, is to consider the full magnetic vector
potential evolution equation in the $xy$-plane:
\begin{align}
\label{eqn:system2.5a}
\A^1_{, t} - u^2 \A^2_{, x} - u^3 \A^3_{, x} + u^2 \A^1_{, y} & =
0, \\
\label{eqn:system2.5b}
\A^2_{, t} + u^1 \A^2_{, x} - u^1 \A^1_{, y} - u^3 \A^3_{, y} & =
0, \\
\label{eqn:system2.5c}
\A^3_{, t} + u^1 \A^3_{, x} + u^2 \A^3_{, y} & = 0.
\end{align}
Approximate solutions to these equations can be computed via the
method outlined in \S\ref{sec:method_potential} and
\S\ref{subsec:mp_limiters}. 
The magnetic field can be obtained via discrete analogs of
\[
 B^1 =  \A^3_{,y}, \quad
 B^2 = - \A^3_{,x}, \quad \text{and} \quad
 B^3 = \A^2_{, x} - \A^1_{, y}.
 \]
Because the numerical update based on the scheme presented
 in \S\ref{sec:method_potential} 
 and \S\ref{subsec:mp_limiters} is not equivalent to the 
 method of \cite{article:Ro04b} in 2.5D, we expect to see 
 some differences in the solution of the magnetic field.
 However, since the two methods are very closely related, we
 expect to obtain  similar results.
 Numerical examples for the 2.5D case comparing the 2D constrained
 transport method of  \cite{article:Ro04b} and the proposed 3D method
 are discussed in \S \ref{sec:smooth2.5d} and \S \ref{sec:shock2.5d}.
 

\section{Numerical experiments}
\label{sec:numex}
We present several numerical experiments in this section.
All the numerical tests are carried out in the
 {\sc clawpack} software package \cite{claw}, which can be
 freely downloaded from the web.
 In particular,  the current work has been incorporated into the   {\sc mhdclaw} extension of {\sc clawpack}, which was originally developed by Rossmanith
 \cite{mhdclaw}. This extension can also be freely downloaded from the web.

\subsection{Test cases in 2.5D}
We begin by presenting numerical results for two problems
in 2.5 dimensions (i.e., the solution depends only on the
independent variables $(x,y,t)$, but all three components of the
velocity and magnetic field vectors are non-zero). The two
examples that are considered are:
\begin{enumerate}
\item Smooth Alfv$\acute{\text{e}}$n wave; and
\item Cloud-shock interaction.
\end{enumerate}
The first example involves an infinitely smooth exact solution
of the ideal MHD system. The second example involves 
the interaction of a high-pressure shock with a
high-density bubble (i.e., a discontinuous example).

\subsubsection{Smooth Alfv\'en wave problem}
\label{sec:smooth2.5d}
We first verify the order of convergence of the constrained transport
method for a smooth circular polarized Alfv\'en wave that propagates
in direction $\n = (\cos \phi,
\sin \phi, 0)^T$ towards the origin. This problem has been considered
by several authors (e.g., \cite{article:Ro04b,article:To00}), and is
a special case ($\theta=0$) of the 3D problem described in detail in
\S \ref{sec:alfven}. Here we take $\phi = \tan^{-1}(0.5)$ and solve
on the domain: 
\begin{equation}
(x,y) \in \left[0,\frac{1}{\cos\phi}\right]\times \left[0,\frac{1}{\sin\phi}\right].
\end{equation}
The solution consists of a sinusoidal wave propagating at constant
speed without changing shape, thus making it a prime candidate
to verify order of accuracy.

In Table \ref{table:alfven2.5}, we show the convergence rates for the
2.5D test problem. Here all three
components of the magnetic field have been updated in the constrained
transport method and the magnetic potential was approximated using a
dimensional split method for (\ref{eqn:system2.5a})--(\ref{eqn:system2.5c}).
The table clearly shows that the proposed method is second-order
accurate in all of the magnetic field components, as well as the magnetic potential
components. The reported errors are of the same magnitude as 
those reported in Table 7.1 of \cite{article:Ro04b} (page 1786).

\begin{table}
\begin{center}
\begin{Large}
  \begin{tabular}{|c||c|c|c|}
    \hline
    {\normalsize {\bf Mesh}} & 
    {\normalsize {\bf $L_{\infty}$ Error in }$B^1$} & 
    {\normalsize {\bf $L_{\infty}$ Error in }$B^2$} & 
    {\normalsize {\bf $L_{\infty}$ Error in }$B^3$}  \\
    \hline \hline
    {\normalsize $64 \times 128$}  &
    {\normalsize $6.778 \times 10^{-4}$}  &  
    {\normalsize $2.393 \times 10^{-3}$}  &  
    {\normalsize $1.284 \times 10^{-2}$}   \\
    \hline
    {\normalsize $128 \times 256$}  &  
    {\normalsize $1.690 \times 10^{-4}$}  &  
    {\normalsize $5.969\times 10^{-4}$}   & 
    {\normalsize $3.203\times 10^{-3}$}   \\
    \hline
    {\normalsize $256 \times 512$}  &  
    {\normalsize $4.221 \times 10^{-5}$}  & 
    {\normalsize $1.492\times 10^{-4}$}  &  
    {\normalsize $8.004\times 10^{-4}$}   \\
    \hline
    {\normalsize $512 \times 1024$}  &  
    {\normalsize $1.057 \times 10^{-5}$}  & 
    {\normalsize $3.741 \times 10^{-5}$}  &  
    {\normalsize $2.000\times 10^{-4}$}   \\
    \hline \hline 
    {\normalsize {\bf Order}} & 
    {\normalsize 2.001} & 
    {\normalsize 2.003} & 
    {\normalsize 2.005} \\
 \hline
  \end{tabular}
  
  \bigskip

  \begin{tabular}{|c||c|c|c|}
    \hline
    {\normalsize {\bf Mesh}} &  
    {\normalsize {\bf $L_{\infty}$ Error in }$\A^1$} & 
    {\normalsize {\bf $L_{\infty}$  Error in }$\A^2$} & 
    {\normalsize {\bf $L_{\infty}$ Error in }$\A^3$} \\
    \hline \hline
    {\normalsize $64 \times 128$}  &  
    {\normalsize $1.302 \times 10^{-2}$} & 
    {\normalsize $1.288\times 10^{-2}$}  &  
    {\normalsize $1.453 \times 10^{-2}$} \\
    \hline
    {\normalsize $128 \times 256$}  &  
    {\normalsize $3.260\times 10^{-3}$}  &
    {\normalsize $3.217\times 10^{-3}$}  &
    {\normalsize $3.633\times 10^{-3}$} \\
\hline
    {\normalsize $256 \times 512$}  &
    {\normalsize $8.213\times 10^{-4}$}  &  
    {\normalsize $8.025\times 10^{-4}$}  &  
    {\normalsize $9.081\times 10^{-4}$} \\
\hline
    {\normalsize $512 \times 1024$}  &    
    {\normalsize $2.089\times 10^{-4}$}  &  
    {\normalsize $1.997\times 10^{-4}$}  &  
    {\normalsize $2.270\times 10^{-4}$} \\
    \hline \hline 
   {\normalsize {\bf Order}} & 
   {\normalsize 1.987} & 
   {\normalsize 2.004} & 
   {\normalsize 2.001} \\
 \hline
  \end{tabular}
    \end{Large}
  \caption{Error tables for the 2.5D Alfv$\acute{\text{e}}$n problem
       at time $t=1.5$ using proposed dimensional split scheme.
       This table shows that all components of the magnetic field and
       all components of the magnetic potential converge at second
       order accuracy.\label{table:alfven2.5}}
\end{center}
\end{table}

\subsubsection{Cloud-shock interaction problem}
\label{sec:shock2.5d}
Next we consider a problem with a strong shock interacting with
a relatively large density jump in the form of a high-density bubble
that is at rest with its background before the shock-interaction.
This problem has also been considered
by several authors (e.g., \cite{article:Ro04b,article:To00}), however,
in previous work a solution was always obtained by treating $B^3$
as a standard conserved variable. Here we  compare such
an approach in the form of the method described in \cite{article:Ro04b},
against the method that is proposed in the current work. In our new method
all three components of the magnetic field are computed from derivatives
of the magnetic vector potential.
The details of the initial conditions are written in \S \ref{sec:cloud-shock}.

\begin{figure}[!t]
   (a)\includegraphics[width=55mm]{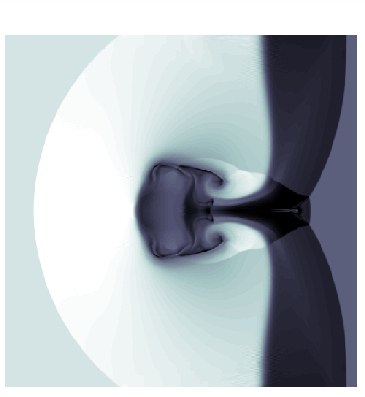}
   (b)\includegraphics[width=55mm]{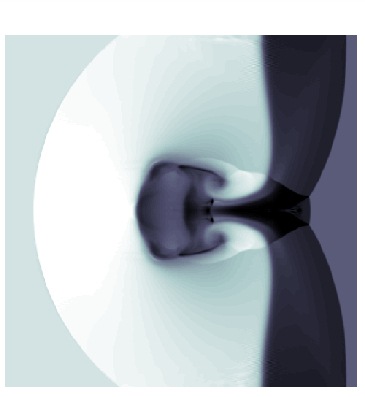}
  \caption{The 2.5D cloud-shock interaction problem. Shown here
  are the out-of-plane magnetic field at time $t=0.06$ 
  as calculated on a  512$\times$512
  mesh using (a) the
  scheme of \cite{article:Ro04b} that only uses $\A^3$, and (b)
 the proposed scheme using the full vector potential $\A$.
\label{figure:cloud-shock2.5}}
\end{figure}

In Figure \ref{figure:cloud-shock2.5} we show numerical results for the
2.5 dimensional problem. For the newly proposed method we take
$\nu = 0.02$ as our artificial diffusion parameter. We compare the scheme from
\cite{article:Ro04b} that only updates $B^1$ and $B^2$ using the
scalar evolution equation for the vector potential (\ref{eqn:mp3})
with the new constrained transport method that updates all components of the
magnetic field.  Although these methods compute $B^3$ in very different ways,
 using very different limiters, these plots show that they nonetheless
 produce similar solutions. This result gives us confidence that the proposed
 method is able to accurately resolve strong shocks despite the somewhat 
 unusual approach for approximating the magnetic field.

\subsection{Test cases in 3D}
We present numerical results for three-dimensional versions of three
classical MHD test problems:
\begin{enumerate}
\item Smooth Alfv\'en wave;
\item Rotated shock tube problem;
\item Orszag-Tang vortex; and
\item Cloud-shock interaction.
\end{enumerate}
 Two-dimensional versions of these
problems have been studied by many authors, see for instance
\cite{article:Ro04b,article:To00}.

 Our test calculations will all be based on the
splitting (\ref{eqn:splitting2a})--(\ref{eqn:splitting2c}). We also carried out several
tests with (\ref{eqn:splitting1a})--(\ref{eqn:splitting1c}); however, we will not report those here.
The two methods produced comparable results.

\subsubsection{Smooth Alfv\'en wave problem}
\label{sec:alfven}
We first verify the order of convergence of the constrained transport
method for a smooth circular polarized Alfv\'en wave that propagates
in direction $\n = (\cos \phi \cos \theta,
\sin \phi \cos \theta, \sin \theta)^T$
towards the origin. Here $\phi$ is an angle with
respect to the $x$-axis in the $xy$-plane and $\theta$ is an angle with respect to the
$xz$-plane. Initial values for the velocity and the magnetic field
are specified in the direction $\n$ as well as the orthonormal
directions $\t = (-\sin \phi, \cos \phi, 0)^T$ and ${\bf r} = (- \cos \phi
\sin \theta, -\sin \phi \sin \theta, \cos \theta)^T$:
\begin{align}
\u (0, \x) &= u^n \, \n + u^t \, \t + u^r \, {\bf r}, \\
\B(0, \x) &= B^n \, \n + B^t \, \t + B^r \, {\bf r},
\end{align}
where
\begin{align}
 u^n &= 0, \quad B^n = 1, \\
 u^t &= B^t = 0.1 \sin(2 \pi \, \n \cdot \x), \\
 u^r &= B^r = 0.1 \cos(2 \pi \, \n \cdot \x).
\end{align}
The initial density and pressure are constant and set to 
\begin{equation}
\rho(0, \x) = 1 \qquad \text{and} \qquad
p(0, \x) = 0.1, 
\end{equation}
respectively.
This choice guarantees that 
the Alfv\'en wave speed is $|v_A|=B_n / \sqrt{\rho} = 1$, which
means that the flow agrees with the initial state whenever the time is an integer
value.
The computational domain is taken to be
\begin{equation}
\Omega = \left[ 0,\frac{1}{\cos \phi \cos \theta} \right] \times \left[ 0, \frac{1}{\sin \phi \cos
  \theta} \right] \times \left[ 0, \frac{1}{\sin \theta} \right],
\end{equation}
with periodic boundary conditions imposed on the conserved
variables in all three coordinate directions: $\left( \rho, \, \rho \u, \, \En, \, \B \right)$.

 Our initial condition for the magnetic potential is
\begin{align}
 \A^1(0, \x) &= 
 	z \left< B^2 \right>  - \frac{1}{20 \pi} \sin \phi
\sin(2 \pi \, \n \cdot \x), \\
 \A^2(0, \x) &= x \left< B^3 \right> + \frac{1}{20 \pi} \cos \phi \sin(2 \pi \, \n
\cdot \x), \\
 \A^3(0, \x) &= y \left< B^1 \right>+ \frac{1}{20 \pi \cos \theta} \cos(2 \pi \, \n
\cdot \x),
\end{align}
where $\left< \B \right>$ denotes the average magnetic field
over the computational domain:
\begin{equation}
\left< \B \right>:=\frac{1}{|\Omega|} \iiint_{\Omega} \B(t,x,y,z) \, dx \, dy \, dz
= \left( \cos \phi \cos \theta, \, \sin \phi \cos \theta, \, \sin \theta \right)^T.
\end{equation}
We note that even though the magnetic field
is time-dependent, its average, $\left< \B \right>$, is time-independent.
Therefore, the magnetic potential consists of a linear (time-independent)
and a periodic (time-dependent) part. Boundary conditions for the 
magnetic potential are handled by applying periodic boundary
conditions on the periodic part and linear extrapolation
for the linear part (linear extrapolation is exact in this case).

In Table \ref{table:alfven3} we show the results of a numerical
convergence study in the three-dimensional
case with $\theta = \phi =
\tan^{-1} (0.5) \approx 26.5651^\circ$. These results confirm that our method is second order accurate.
\begin{table}
\begin{center}
\begin{Large}
  \begin{tabular}{|c||c|c|c|}
    \hline
    {\normalsize {\bf Mesh}} & 
    {\normalsize {\bf $L_{\infty}$ Error in }$B^1$} & 
    {\normalsize {\bf $L_{\infty}$ Error in }$B^2$} & 
    {\normalsize {\bf $L_{\infty}$ Error in }$B^3$}  \\
    \hline \hline
    {\normalsize $16 \times 32 \times 32$}  &   
    {\normalsize $ 1.022 \times 10^{-2}$}  &  
    {\normalsize $2.787 \times 10^{-2}$}  &  
    {\normalsize $2.382 \times 10^{-2}$}   \\
    \hline
    {\normalsize $32 \times 64 \times 64$}  &  
    {\normalsize $2.577 \times 10^{-3}$}  &  
    {\normalsize $7.075\times 10^{-3}$}   & 
    {\normalsize $6.101\times 10^{-3}$}   \\
    \hline
    {\normalsize $64 \times 128 \times 128$}  &  
    {\normalsize $6.487 \times 10^{-4}$}   & 
    {\normalsize $1.782\times 10^{-3}$}  &  
    {\normalsize $1.549\times 10^{-3}$}   \\
    \hline \hline 
 {\normalsize {\bf Order}} & 
 {\normalsize 1.990} & 
 {\normalsize 1.989} & 
 {\normalsize 1.978} \\
 \hline
  \end{tabular}
  
  \bigskip
  
  \begin{tabular}{|c||c|c|c|}
    \hline
    {\normalsize {\bf Mesh}} &  
    {\normalsize {\bf $L_{\infty}$ Error in }$\A^1$} & 
    {\normalsize {\bf $L_{\infty}$  Error in }$\A^2$} & 
    {\normalsize {\bf $L_{\infty}$ Error in }$\A^3$} \\
    \hline \hline
    {\normalsize $16 \times 32 \times 32$}  &     
    {\normalsize $1.902 \times 10^{-1}$} & 
    {\normalsize $1.460\times 10^{-1}$}  &  
    {\normalsize $1.187 \times 10^{-1}$} \\
    \hline
    {\normalsize $32 \times 64 \times 64$}  &  
    {\normalsize $4.816\times 10^{-2}$}  &  
    {\normalsize $3.747\times 10^{-2}$}  &  
    {\normalsize $2.995\times 10^{-2}$} \\
    \hline
    {\normalsize $64 \times 128 \times 128$}  &    
    {\normalsize $1.220\times 10^{-2}$}  &  
    {\normalsize $9.528\times 10^{-3}$}  &  
    {\normalsize $7.564\times 10^{-3}$} \\
    \hline \hline 
    {\normalsize {\bf Order}} & 
    {\normalsize 1.981} & 
    {\normalsize 1.975} & 
    {\normalsize 1.985} \\
 \hline
  \end{tabular}
  \end{Large}
  \caption{\label{table:alfven3}
Error tables for the 3D Alfv$\acute{\text{e}}$n
    problem at time $t=1$.
  	This table shows that all components of the magnetic field and
	all components of the magnetic potential converge at second
	order accuracy. The approximate order of accuracy is computed
	by taking the $\log_2$ of the ratio of the error on the
	$32\times64\times64$ mesh to the error on the 
	$64\times128\times128$ mesh.}
\end{center}
\end{table}

\subsection{Rotated shock tube problem}
We now consider the shock tube problem described in 
\cite{article:MiTz10,article:StGa09}. In this test case a one-dimensional shock tube problem is computed in a three-dimensional domain. In order to define the problem,
we consider a coordinate transformation of the form 
\begin{gather}
\label{eqn:rp_transformation}
  \begin{bmatrix}
    \xi \\ \eta \\ \zeta 
  \end{bmatrix} =
\begin{bmatrix}
  \cos(\alpha)  \cos(\beta)  & \cos(\alpha) \sin(\beta) & \sin(\alpha) \\
   - \sin(\beta) & \cos(\beta) & 0 \\
   -\sin(\alpha)  \cos(\beta) &-\sin(\alpha)  \sin(\beta)  & \cos(\alpha)
\end{bmatrix}
  \begin{bmatrix} x \\ y \\ z 
  \end{bmatrix} ,  
\end{gather}
which maps the Cartesian coordinates $(x,y,z)^T$ to the rotated
coordinate system  $(\xi,\eta,\zeta)^T$.
Vectors transform as follows
\begin{gather}
\label{eqn:rp_transform2}
  \begin{bmatrix}
    B^x \\ B^y \\ B^z 
  \end{bmatrix} =
\begin{bmatrix}
  \cos(\alpha)  \cos(\beta)  & -\sin(\beta) & -\sin(\alpha)  \cos(\beta) \\
   \cos(\alpha)  \sin(\beta) & \cos(\beta) & -\sin(\alpha)  \sin(\beta) \\
  \sin(\alpha)  &  0  & \cos(\alpha)
\end{bmatrix}
  \begin{bmatrix} B^{\xi} \\ B^{\eta} \\ B^{\zeta} 
  \end{bmatrix}.
\end{gather}

The Riemann initial data is taken to be
\begin{equation}
\begin{split}
&\left(\rho,\, u^{\xi},\, u^{\eta}, \, u^{\zeta}, \, p, \, B^{\xi}, \, B^{\eta}, \, B^{\zeta} \right)(0, \x) \\
&=\left\{ \begin{array}{l c c}
\Bigl(1.08, \, 1.2, \, 0.01, \, 0.5, \, 0.95, \, \frac{2}{\sqrt{4\pi}}, \, \frac{3.6}{\sqrt{4\pi}}, \, 
\frac{2}{\sqrt{4\pi}} \Bigr)
& \text{if} & \xi < 0,\\
\Bigl( 1,\, 0,\, 0,\, 0,\, 1,\, \frac{2}{\sqrt{4\pi}}, \, \frac{4}{\sqrt{4\pi}}, \, \frac{2}{\sqrt{4 \pi}} \Bigr) & \text{if} & \xi \ge 0,
\end{array} \right.
\end{split}
\end{equation}
where the velocity vector and the magnetic field are given in the
rotated coordinate frame.
As the initial condition for the magnetic potential (in the rotated coordinate
frame) we use
\begin{equation}
\left({\A}^\xi, \, {\A}^\eta, \, {\A}^\zeta \right)(0,\x)
=\left(0, \, \xi B^{\zeta}, \, \eta B^{\xi} - \xi B^{\eta}\right).
\end{equation}

The solution of the Riemann problem at later times remains a function
of $\xi$. This fact is used to extend the computed values to ghost
cell values. The computational domain is 
$\Omega = [-0.75, 0.75] \times [0, 0.015625] \times [0, 0.015625]$,
which is discretized using $768 \times 8 \times 8$ grid cells. By using
$\alpha = \tan^{-1}(0.5)$ and $\beta = \tan^{-1}(0.25 \cos \alpha)$, we
obtain simple formulas for the extension of the numerical solution to
the ghost cells.To define ghost cell values in the Cartesian $y$-direction, we
can for instance use the relation $q(i,j,k) = q(i+1, j-2, k)$. In the
$z$-direction we can use $q(i,j,k) = q(i+1,j,k-4)$. The ghost cell
values in the $x$-direction were computed via extrapolation.
We present numerical results at time $t=0.2$ as scatter plots of the
three-dimensional solution plotted  as a function of $\xi$. 

\begin{figure}[htb]  
(a)\includegraphics[width=0.5\textwidth]{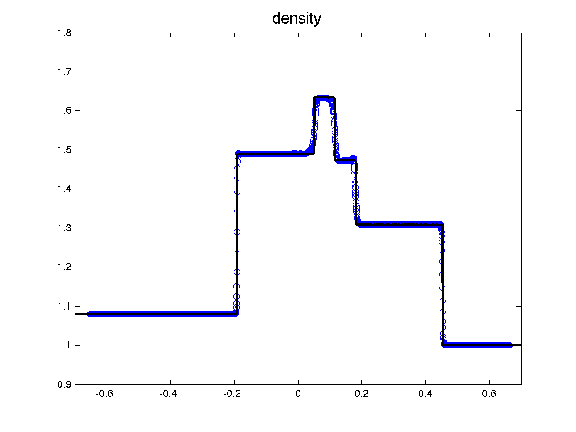}\hfil
(b)\includegraphics[width=0.5\textwidth]{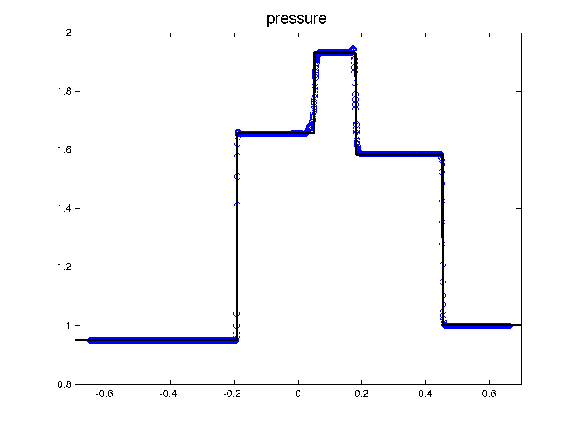}

(c)\includegraphics[width=0.5\textwidth]{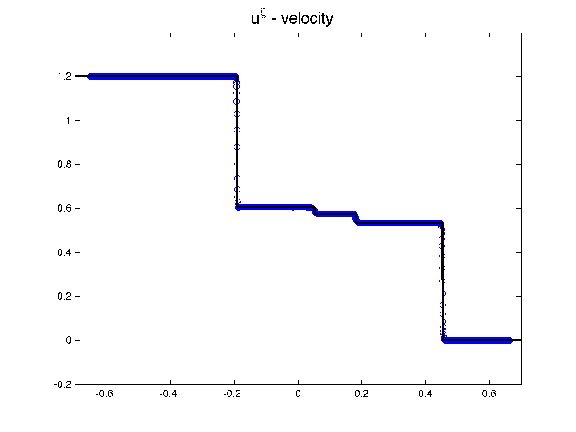}\hfil
(d)\includegraphics[width=0.5\textwidth]{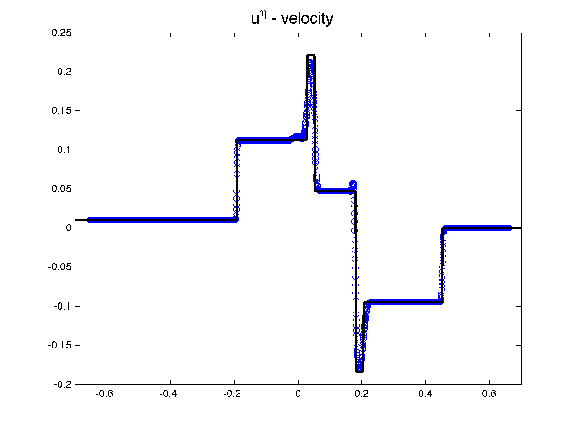}

(e)\includegraphics[width=0.5\textwidth]{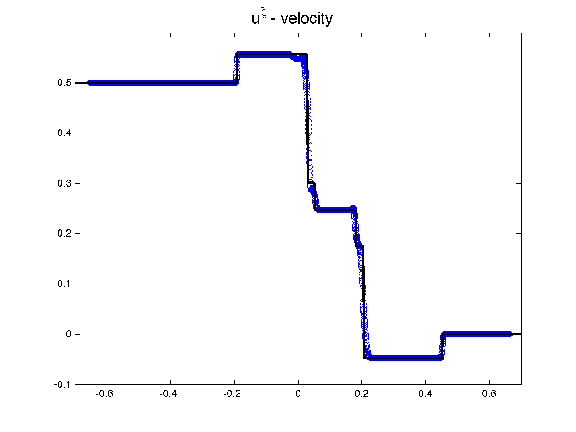}\hfil
(f)\includegraphics[width=0.5\textwidth]{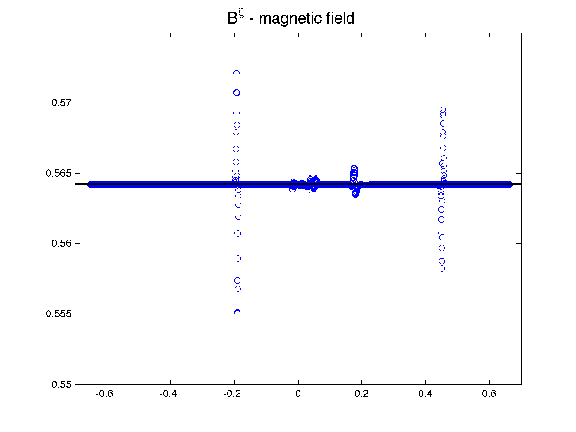}

(g)\includegraphics[width=0.5\textwidth]{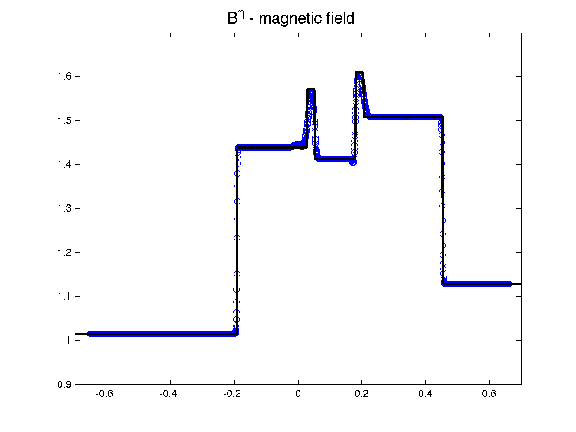}\hfil
(h)\includegraphics[width=0.5\textwidth]{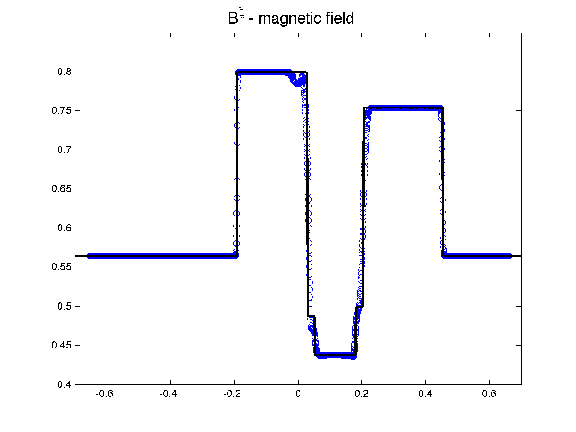}
  \caption{\label{fig:rp} The rotated Riemann problem. Scatterplots of various components of  solution at time $t=0.2$
    using the constrained transport algorithm. The solid line in each panel is a highly-resolved
    1D simulation.}
\end{figure}

As already documented by
other authors, we also observe some oscillations in particular in the $B^{\xi}$ 
component of the magnetic field. These oscillations
do not appear if the shock tube initial data are aligned with
the mesh.   
In Figure \ref{fig:rp} we show approximations of the solution at time
$t=0.2$. These plots are scatterplots, where the output of the
three-dimensional computation is plotted as a function of $\xi$.
The solid line in these plots is obtained by computing solutions of
the one-dimensional Riemann problem on a very fine mesh (using $10^4$
grid cells).
The computation was performed using the minmod limiter and the diffusive
limiter with $\nu=0.05$. As expected, with a less diffusive limiter the
amplitude of the oscillations becomes stronger, 
by increasing $\nu$, the amplitude of the oscillations becomes weaker. 

%
%
%


\subsubsection{Orszag-Tang vortex}
Our next test problem is the Orszag-Tang vortex problem. This is a
standard test case for the two-dimensional MHD equations. Here we
add a  small perturbation to the initial velocity field that depends
on the vertical direction.

The initial condition is
\begin{align}
\rho(0, \x) &= \gamma^2, \quad p(0, \x) = \gamma, \\
\u(0, \x) & = \Bigl(- (1+\varepsilon \sin z)  \sin y, \,  (1+ \varepsilon
\sin z) \sin x, \, \varepsilon \sin(z) \Bigr)^T,\\
\B(0, \x) & =\Bigl(- \sin y, \, \sin(2x), \, 0 \Bigr)^T.
\end{align}  
Here we used $\varepsilon = 0.2$. 
The initial condition for the magnetic potential is 
\begin{equation}
\Av(0, \x) = \Bigl( 0, \, 0,  \, \cos y + \cos(2 x) \Bigr)^T.
\end{equation}
The computational domain is a
cube with side length $2 \pi$. Periodicity is imposed in all three
directions.

\begin{figure}[htb]  
 \includegraphics[width=1\textwidth]{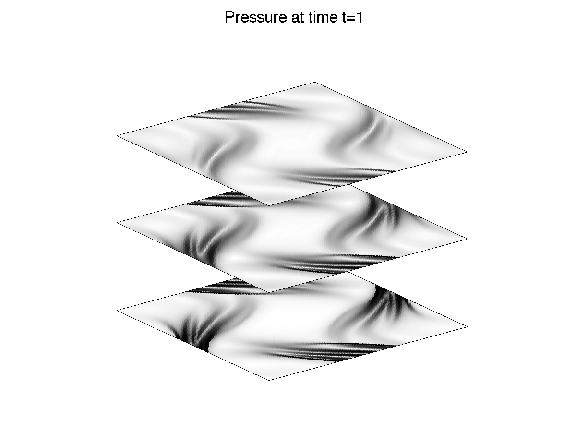}
  \caption{\label{fig:ot_test2_a} The 3D Orszag-Tang problem.
  Shown in this figure are schlieren slices of pressure at time $t=1$ 
    and at various $z$ values
    ($z=\pi/2, \, \pi, \, 3\pi/2$). This solution was computed on a
    mesh with $150 \times 150 \times 150$ grid cells.}
\end{figure}

\begin{figure}[htb]  
\includegraphics[width=1\textwidth]{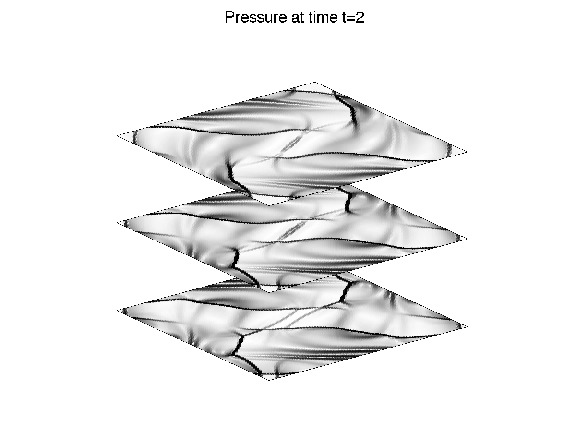}
  \caption{\label{fig:ot_test2_b} 
    The 3D Orszag-Tang problem.
  Shown in this figure are schlieren slices of pressure at time $t=2$ 
    and at various $z$ values
    ($z=\pi/2, \, \pi, \, 3\pi/2$). This solution was computed on a
    mesh with $150 \times 150 \times 150$ grid cells.}
\end{figure}

\begin{figure}[htb]  
\includegraphics[width=1\textwidth]{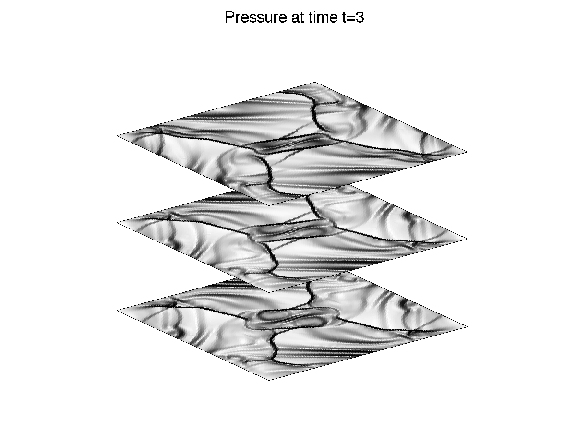}
  \caption{\label{fig:ot_test2_c} 
  The 3D Orszag-Tang problem.
  Shown in this figure are schlieren slices of pressure at time $t=3$ 
    and at various $z$ values
    ($z=\pi/2, \, \pi, \, 3\pi/2$). This solution was computed on a
    mesh with $150 \times 150 \times 150$ grid cells.}
\end{figure}

\begin{figure}[htb]  
\includegraphics[width=1\textwidth]{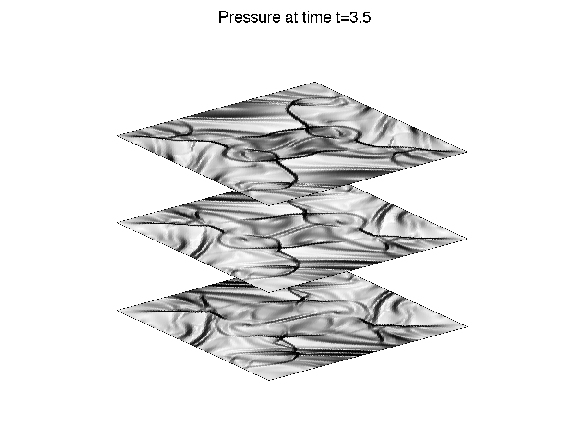}
  \caption{\label{fig:ot_test2_d}
  The 3D Orszag-Tang problem.
  Shown in this figure are schlieren slices of pressure at time $t=3.5$ 
    and at various $z$ values
    ($z=\pi/2, \, \pi, \, 3\pi/2$). This solution was computed on a
    mesh with $150 \times 150 \times 150$ grid cells.}
\end{figure}

\begin{figure}[htb]  
 (a)\includegraphics[width=0.9\textwidth]{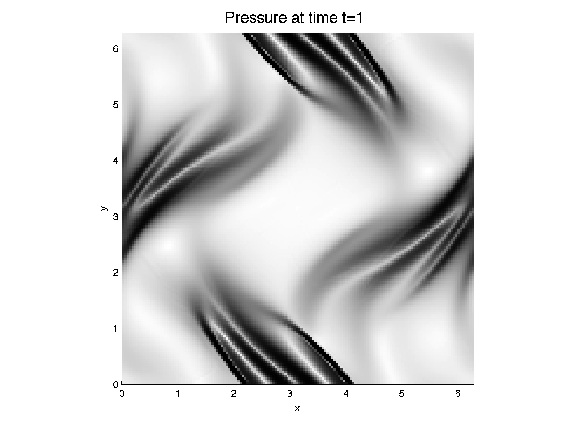}

(b)\includegraphics[width=0.9\textwidth]{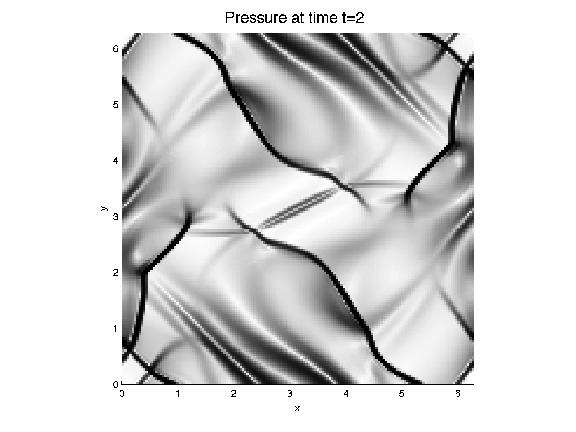}
  \caption{\label{fig:ot_test2_xy_a} 
  The 3D Orszag-Tang problem. Schlieren plots of pressure 
  at time (a) $t=1$ and (b) $t= 2$ computed on a
    mesh with $150 \times 150 \times 150$ grid cells. Slices
    of the solution at $z=\pi$ in the $xy$-plane are shown.}
\end{figure}

\begin{figure}[htb]  
(a)\includegraphics[width=0.9\textwidth]{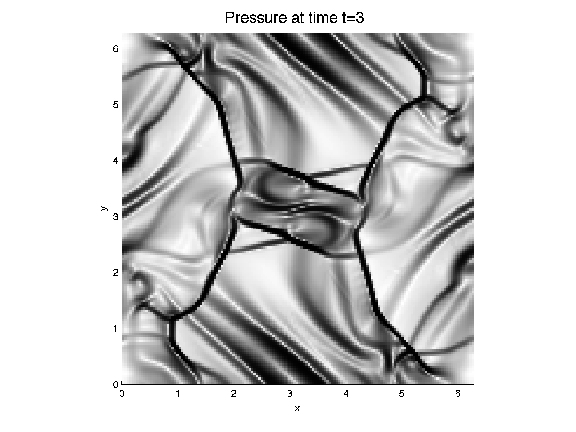}

(b)\includegraphics[width=0.9\textwidth]{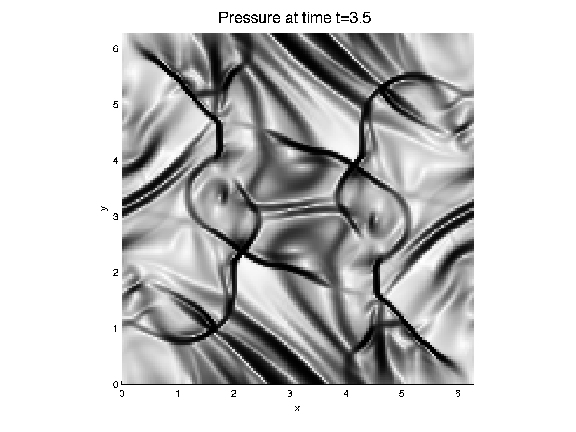}
  \caption{\label{fig:ot_test2_xy_b} The 3D Orszag-Tang vortex problem. Shown
  in these panels are schlieren plots of pressure in the $xy$-plane 
  at (a) time $t=3$ and $z=\pi$ and
  (b) $t=3.5$ and $z=\pi$. These solutions were computed on a
    mesh with $150 \times 150 \times 150$ grid cells.}
\end{figure}

In Figures \ref{fig:ot_test2_a}-\ref{fig:ot_test2_d} we show a
sequence of schlieren plots of the pressure
for several slices of data in the $xy$-plane for $z=\pi/2, \, \pi, \, 3\pi/2$. 
Those computations have been
performed with the high-resolution constrained transport method with 
monotonized central limiter and the diffusive limiter using $\nu =
0.05$.
In Figures \ref{fig:ot_test2_xy_a}, \ref{fig:ot_test2_xy_b} we show
schlieren plots of pressure  in the $xy$-plane for $z=\pi$ at
different times. 

\subsubsection{Cloud-shock interaction problem}\label{sec:cloud-shock}
Finally we consider the cloud-shock interaction problem. 
The initial conditions consist of a shock that is initially located at
$x=0.05$:
\begin{equation}
\begin{split}
&\left(\rho, \, u^1, \, u^2, \, u^3,  \, p,  \, B^1, \, B^2, \, B^3 \right)(0, \x) \\
&=\biggl\{ \begin{array}{l c c}
\left(3.86859,  \, 11.2536,  \, 0,  \, 0,  \, 167.345,  \, 0,  \, 2.1826182, 
 \, -2.1826182\right)
& \text{if} &x < 0.05,\\
\left( 1,  \, 0,  \, 0,  \, 0,  \, 1,  \, 0,  \, 0.56418958,  \, 0.56418958\right) & \text{if} & x\ge 0.05,
\end{array}
\end{split}
\end{equation}
and a spherical cloud of density $\rho = 10$ with radius $r=0.15$ and
centered at $(0.25,0.5,0.5)$. The cloud is in hydrostatic
equilibrium with the fluid to the right of the shock.
The initial condition for the magnetic potential is
\begin{equation}
\Av (0, \x) = \left\{ \begin{array}{l c c}
\left( 2.1826182  \,  y,  \, 0,  \,  -2.1826182 \,  (x-0.05) \right)^T & \text{if} & x < 0.05,
\\
\left( -0.56418958 \,  y,  \, 0,  \, 0.56418956 \, (x-0.05) \right)^T & \text{if} & x \ge 0.05.
\end{array}\right.
\end{equation}

This test problem can be computed on the unit cube with inflow
boundary conditions on the left side and outflow conditions on all
other sides. Instead we make use of the symmetry of the problem and
compute the problem only in a quarter of the full domain, i.e.\ in 
$[0, 1] \times [0.5, 1] \times [0.5, 1]$ and impose reflecting
boundary conditions at the lower boundary in the $y$ and $z$-directions.
For the conserved quantities, the reflecting boundary condition is
implemented in the standard way, i.e., by copying the values of the
conserved quantities from the first grid cells of the flow domain to
the ghost cells. The normal momentum component in the ghost cells
is negated. The components of the magnetic potential are linearly
extrapolated to the ghost cells.
 
\begin{figure}[htb]
\begin{center}
(a)\includegraphics[width=0.9\textwidth]{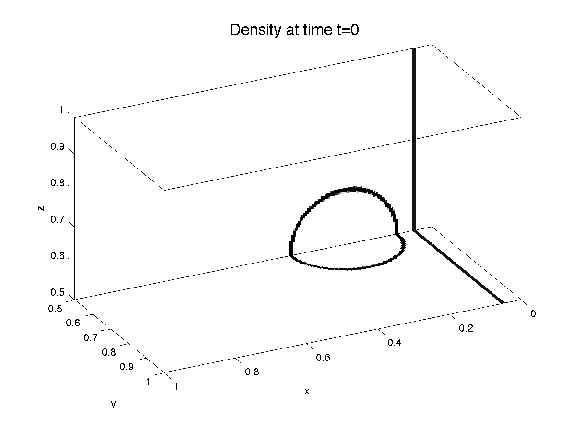}
 
(b)\includegraphics[width=0.9\textwidth]{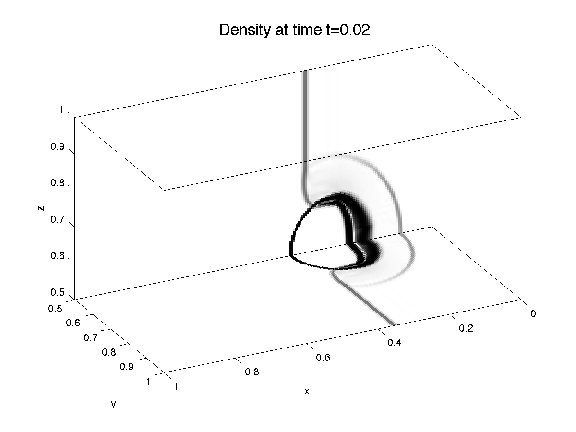}
  \caption{The 3D cloud-shock interaction problem.
  Schlieren plots of density for the problem using a mesh with
    $200\times100\times100$ grid cells at time (a) $t=0$ and (b) $t=0.02$. 
    Shown here is the solution in
    two selected orthogonal planes. 
\label{figure:cloud-shock3_density1}}
\end{center}
\end{figure}
\begin{figure}[htb]
\begin{center}
(a)\includegraphics[width=0.9\textwidth]{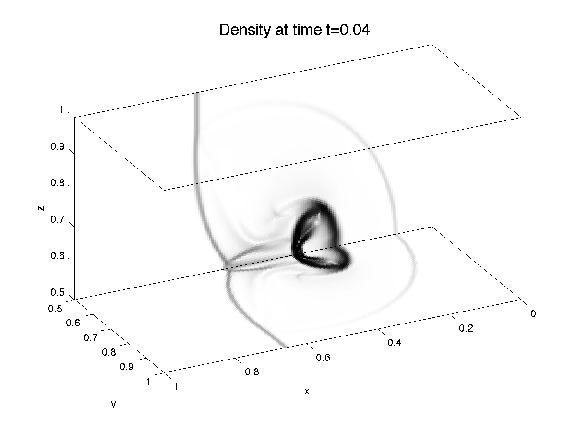}
 
(b)\includegraphics[width=0.9\textwidth]{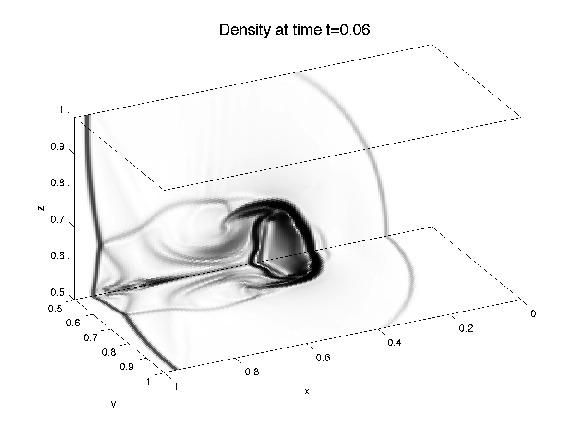}
  \caption{The 3D cloud-shock interaction problem. 
  Schlieren plots of density using a mesh with
    $200\times100\times100$ grid cells at time (a) $t=0.04$ and (b) $t=0.06$. 
    Shown here is the solution in
    two selected orthogonal planes. 
\label{figure:cloud-shock3_density2}}
\end{center}
\end{figure}
\begin{figure}[htb]
\begin{center}
(a)\includegraphics[width=0.9\textwidth]{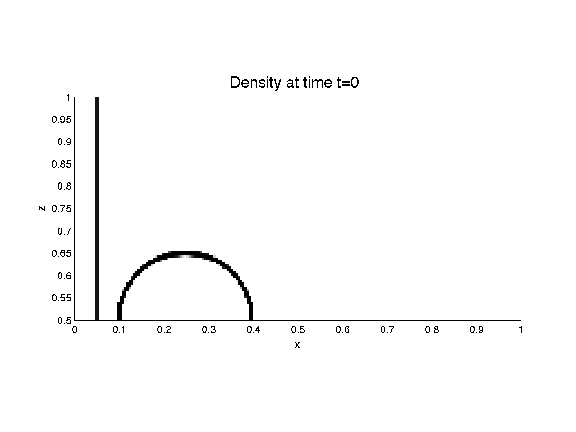}

(b)\includegraphics[width=0.9\textwidth]{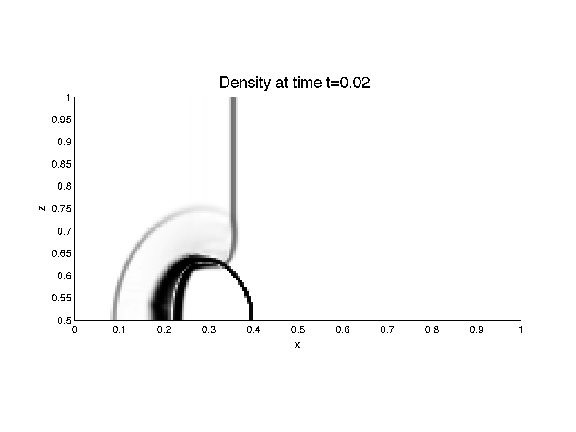}
\caption{\label{figure:cloud-shock_xz_a}
The 3D cloud-shock interaction problem. 
Sequence of schlieren plots of density in the $xz$-plane
  for $y=0.5$ at time (a) $t=0$ and (b) $t=0.02$.}
\end{center} 
\end{figure}
\begin{figure}[htb]
\begin{center}
(a)\includegraphics[width=0.9\textwidth]{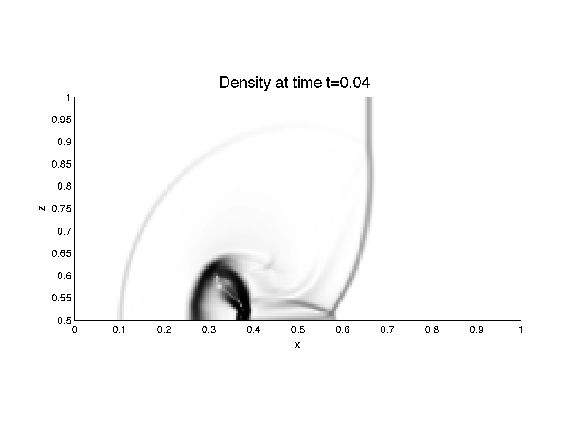}

(b)\includegraphics[width=0.9\textwidth]{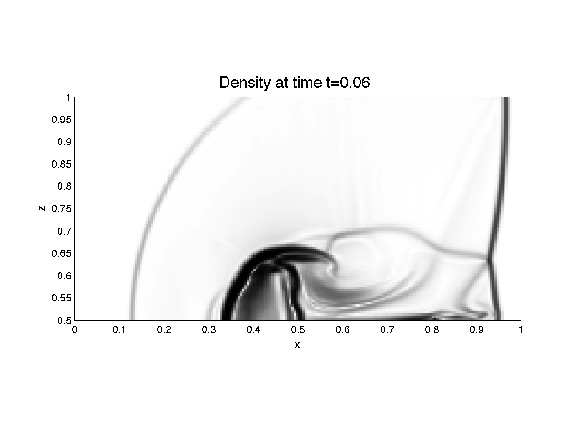}
\caption{\label{figure:cloud-shock_xz_b}
The 3D cloud-shock interaction problem. 
Sequence of schlieren plots of density in the $xz$-plane
  for $y=0.5$ at time (a) $t=0.04$ and (b) $t=0.06$.}
\end{center}
\end{figure}

In Figures \ref{figure:cloud-shock3_density1} and \ref{figure:cloud-shock3_density2}, we show  a sequence of 
schlieren plots of the
density in the $xy$-plane for $z=0.5$ and the $xz$-plane for $y=0.5$. 
The three-dimensional radial symmetric solution structure compares well with
previously shown two-dimensional computations.
In Figures \ref{figure:cloud-shock_xz_a} and
\ref{figure:cloud-shock_xz_b} 
we show the schlieren plots of
density  in the $xz$-plane. 
Here the diffusive limiter described in \S \ref{sec:method_potential} was used
with $\nu = 0.02$.

\section{Conclusions}
We have presented a new constrained transport method for the
three-dimen\-sional MHD equations. 
Depending on the gauge condition, we discussed different
possible evolution equations for the magnetic potential. 
All of these gauge conditions implicate some difficulties for
the discretization. The condition used here leads to a weakly
hyperbolic system for the evolution of the magnetic potential. We
discretized this system with a splitting approach.
For the MHD equations we used a wave propagation method.
Several numerical tests confirm the robustness and
accuracy of the resulting constrained transport scheme. 

Our method is fully explicit, as well as fully unstaggered, and
therefore well-suited for adaptive mesh refinement and parallelization.
These generalizations will be the focus of future work.

\bigskip

\bigskip

\noindent
{\bf Acknowledgements.} The authors would like to thank Prof. Dr. Rainer Grauer
from the Ruhr-Universit\"at-Bochum for several useful discussions, as
well as the anonymous reviewers for their valuable comments.
This work was supported in part by NSF grants DMS-0711885 and 
DMS-1016202 and by the DFG through FOR1048.

\end{document}